\documentclass[a4paper,11pt,leqno]{article}
\usepackage{theorem}
\usepackage{latexsym}
\usepackage{amsmath}
\usepackage{delarray}
\usepackage{graphicx}
\usepackage{delarray}
        \usepackage[latin1]{inputenc}
        \usepackage{color}
        \usepackage{array}
        \usepackage{longtable}
        \usepackage{calc}
        \usepackage{multirow}
        \usepackage{hhline}
        \usepackage{ifthen}

\makeatletter
\@addtoreset{equation}{section}
\makeatother

\newtheorem{cor}{Corollary}[section]
\newtheorem{thm}{Theorem}[section]
\newtheorem{lem}{Lemma}[section]

\theorembodyfont{\rmfamily}
\newtheorem{pft}{Proof of Theorem}[section]
\newtheorem{pfl}{Proof of Lemma}[section]

\def\sgn{\mathop{\rm sgn}}

\let\margin\marginpar
\newcommand\myMargin[1]{\margin{
      \raggedright\tiny  \setlength{\baselineskip}{0pt} #1}}
\renewcommand{\marginpar}[1]{\myMargin{#1}}

\title{Fisher Information Matrix of General Stable Distributions Close to the Normal Distribution}
\author{By Muneya MATSUI\footnote{mmatsui@grad.e.u-tokyo.ac.jp} \\ \\
  Graduate School of Economics, University of Tokyo, \\
 Hongo 7-3-1, Bunkyo-Ku, Tokyo 113-0033, Japan }
\date{}

\begin{document}
\maketitle

\begin{abstract}
We investigate behavior of the Fisher
information matrix of general stable distributions. 
DuMouchel (1975, 1983) proved that the Fisher information
 $I_{\alpha\alpha}$ of characteristic exponent $\alpha$
 diverges to infinity as $\alpha$ approaches $2$.    
Nagaev  and  Shkol'nik (1988) made more detailed analysis of
$I_{\alpha\alpha}$ and derived asymptotic behavior of
 $I_{\alpha\alpha}$ diverging to infinity as $\alpha$ approaches $2$
in the symmetric case.
Extending their work in this paper we have obtained behavior of the Fisher
 information matrix of general stable distributions as $\alpha$ approaches $2$ by detailed
study of behavior of the corresponding density and its score functions. 
We clarify the limiting values of the $4\times4$ Fisher Information
 matrix with respect to the location $\mu$, the scale $\sigma$, the 
 characteristic exponent $\alpha$ and the skewness parameter $\beta$.
\end{abstract}
\begin{small}
\begin{center}
{\bf Keywords} \\
General stable distributions, information matrix, normal distribution,
score functions. 
\end{center}
\end{small}
\section{Notations and preliminary results}
\label{sed:preliminaris}
The family of stable distributions has enjoyed great interest of
researchers in many fields e.g., mathematics, physics, cosmology and
even economics. These applications are summarized in Uchaikin and Zolotarev (1999).
In statistical inference estimation of stable parameters has been of
great interest. In recent years, maximum likelihood estimation of stable distributions
has become feasible (see Brorsen ans Yang (1990), Nolan (2001) or Matsui
and Takemura (2004)). Even in time series models like GARCH
using general stable distributions, maximum likelihood estimation is possible owing to recent
development of global algorithm (see Liu and Brorsen (1995)). Since the
Fisher information matrix gives useful criteria for the accuracy of
estimation, it is indispensable to analyze that of general stable
distributions. Near Gaussian
distribution ($\alpha=2$), the information of $\alpha$,
$I_{\alpha\alpha}$ diverges to 
$\infty$ and asymptotic behavior of $I_{\alpha\alpha}$ as
$\alpha \uparrow 2$ is of great interest. Nagaev and Shkol'nik (1988) have
solved this problem excellently for symmetric stable distributions. 
However for general stable distributions Nagaev and Shkol'nik (1988) stated 
``We note the problems under study are as yet unresolved for non-symmetric stable distributions.''
In this paper we investigate and obtain
asymptotic behavior of $I_{\theta\theta},\
\theta=\mu,\sigma,\alpha,\beta$, as $\alpha \uparrow 2$ under general
stable distributions.

Let
\begin{equation*}
\label{eq:characteristic-function}
\Phi(t)=\Phi(t;\mu,\sigma,\alpha,\beta)=
\exp\left(-|\sigma t|^{\alpha}\left\{1 + i\beta(\sgn
t)\tan\left(\frac{\pi\alpha}{2}\right)(|\sigma t|^{1-\alpha}-1)\right\} + i\mu t\right)
\end{equation*}
denote the characteristic function of general stable distribution
($\alpha \neq 1$)
with parameters
\begin{equation*}
\theta=(\theta_1, \theta_2, \theta_3,\theta_4)=(\mu,\sigma,\alpha,\beta),
\end{equation*}
\begin{equation*}
\mu\in \mathbf{R},\ \sigma >0,\ 0<\alpha \le2,\ |\beta|\le 1.
\end{equation*}
For the standard case
$(\mu,\sigma)=(0,1)$ we simply write the characteristic function as
\begin{equation}
\label{eq:characteristic-function-standard}
\Phi(t;\alpha,\beta)=\exp\left(-|t|^{\alpha}\left\{1 + i\beta(\sgn
t)\tan\left(\frac{\pi\alpha}{2}\right)(|t|^{1-\alpha}-1)\right\} \right).
\end{equation}
This is Zolotarev's (M) parameterization (see p.11 of Zolotarev (1986)).
The corresponding density is written as
$f(x;\mu,\sigma,\alpha,\beta)$  and $f(x;\alpha,\beta)$ in the standard
case. Then
$$
f(x;\mu,\sigma,\alpha,\beta) =
\frac{1}{\sigma}f\left(\frac{x-\mu}{\sigma};\alpha,\beta\right).
$$
We also write the density of $N(0,2)$ as
$$
f(x;2)=\frac{1}{2\sqrt{\pi}}\exp\left(-\frac{x^2}{4}\right).
$$
For the rest of this paper, without loss of generality we consider the Fisher information 
only at the standard case $(\mu,\sigma)=(0,1)$. 

The first derivative of $f(x;\alpha,\beta)$ with respect to $x$ 
is denoted by $f'(x;\alpha,\beta)$ and the partial derivative with
respect to $\theta_i$, i.e., 
\begin{equation*}
\frac{\partial f(x;\mu,\sigma,\alpha,\beta)}{\partial \theta_i}\bigg|_{(\mu,\sigma)=(0,1)}
\end{equation*} 
is denoted by $f_{\theta_i}(x;\alpha,\beta)$.
Note that
\begin{equation}
\label{eq:f-mu}
f_\mu(x;\alpha,\beta)=-f'(x;\alpha,\beta),
\end{equation}
\begin{equation}
\label{eq:f-sigma}
f_\sigma(x;\alpha,\beta)= -f(x;\alpha,\beta) - x f'(x;\alpha,\beta).
\end{equation}
The components of Fisher information matrix $I$ are defined as follows.
\begin{equation*}
\label{eq:information}
   I_{ij}=I_{\theta_i \theta_j}=\int_{-\infty}^{\infty} \frac{\partial f}{\partial \theta_i}
\frac{\partial f}{\partial \theta_j} \frac{1}{f} dx.
\end{equation*}

Let 
\begin{equation*}
\zeta=-\beta\tan\left(\frac{\pi\alpha}{2}\right),
\end{equation*}
\begin{equation*}
\label{varrho}
\varrho=\frac{2}{\pi\alpha}\arctan\left(\beta\tan\left(\frac{\pi\alpha}{2}\right)\right).
\end{equation*}
{}From Theorem 1 of Nolan (1997), which is a modified version of the formula (2.2.18) of 
Zolotarev (1986), 
the density $f(x;\alpha,\beta)$ for the case 
$\alpha\neq 1$ and $x\neq\zeta$  is written as
\begin{equation}
\label{eq:dense}
f(x;\alpha,\beta)
 = \frac{\alpha |x-\zeta|^{1/(\alpha-1)}}{2|\alpha-1|}\int^{1}_{-\varrho^{\ast}} 
     A(\varphi;\alpha,\beta)\exp\left(-|x-\zeta|^{\alpha/(\alpha-1)}A(\varphi;\alpha,\beta)\right) d\varphi, 
\end{equation}
where  $\varrho^{\ast}=\varrho \sgn(x-\zeta)$ and
\begin{equation}
A(\varphi;\alpha,\beta)=\left(\cos\left(\frac{\pi}{2}\alpha\varrho\right)\right)^{\frac{1}{\alpha-1}}
 \left(\frac{\cos\left(\frac{\pi}{2}\varphi\right)}{\sin \left\{\frac{\pi}{2}
\alpha (\varphi+\varrho)\right\}
 }\right)^{\frac{\alpha}{\alpha-1}}\frac{\cos\left[\frac{\pi}{2}\{\alpha\varrho+(\alpha-1)\varphi\}\right]}
 {\cos(\frac{\pi}{2}\varphi)} 
\label{integrand-density-A1}
\end{equation}
is a positive function.
Here we should note the relation
$f(x;\alpha,\beta)=f(-x;\alpha,-\beta)$ for $x-\zeta<0$, which will be
used many times.
As $\alpha \uparrow 2$, $\zeta$ and $\varrho$ converge to 0 and at
$\alpha=2$ we obtain an unusual representation of Gaussian distribution $N(0,2)$, 
\begin{equation*}
 f(x;2)=x \int_0^1 1/\left(2\sin\left(\frac{\pi}{2}\varphi\right)\right)^2
\exp\left(-x^2/\left(2\sin\left(\frac{\pi}{2}\varphi\right)\right)^2 \right)\  d\varphi,\quad x>0.
\end{equation*} 

The remainder of this paper consists of three sections. In Section 2 we derive behavior of
$f'(x;\alpha,\beta)$, $f_\mu(x;\alpha,\beta)$ and $f_{\sigma}(x;\alpha)$
as $\alpha \uparrow 2$ by careful analysis of the formula
(\ref{eq:dense}). In Section 3 behavior of $f_\alpha(x;\alpha,\beta)$ and
$f_\beta(x;\alpha,\beta)$ are obtained by utilizing the inversion
formula. The information matrix of general stable distributions is given in Section 4.

\section{Density and derivatives of general stable distributions close
 to the normal distribution}
In this Section the density and derivatives of general stable
distributions as $\alpha \uparrow 2$ and $x \to \infty$ are obtained.
In the following we write
$$
\Delta=2-\alpha.
$$
As $\Delta \to 0$ we can show 
\begin{equation}
\label{normal}
f(x;\alpha,\beta)=f(x;2)+O(\Delta)
\end{equation}
 uniformly in $x$
by the finiteness of
$f_{\alpha}(x;\alpha,\beta)$, which can be easily verified by 
the formula (\ref{eq:differential-alpha-density-inversion-formula-half})
in the proof of Lemma \ref{lem:f_alpha(x;alpha,beta)}. However if we
consider the ratio $f(x;2)/f(x;\alpha,\beta)$ rather than the difference
$f(x;2)-f(x;\alpha,\beta)$ of the densities as $x \to \infty$,
$f(x,2)$ is much smaller than $f(x;\alpha,\beta)$ and we can not obtain an
accurate approximation of $f(x,\alpha,\beta)$ without closer investigation of
the remainder $O(\Delta)$. In the case of $x \to \infty$ we can utilize 
the asymptotic expansion
\begin{equation}
\label{expansion}
f(x;\alpha,\beta)=\Delta(1+\beta)(x-\zeta)^{\Delta-3}+o(\Delta x^{-4}) 
\end{equation}
from the expansion (2.5.4) on p.94 of Zolotarev (1986), Bergstr\"om (1953)
or Section XVII.6 of Feller (1971). 
(2.5.4) of Zolotarev (1986) is in terms of (B) representation with
characteristic function defined on p.12 of Zolotarev (1986). Hence as $x\to\infty$
\begin{equation*}
f_B(x;\alpha,\beta)=\Delta(1+\beta_B)x^{\Delta-3}+o(\Delta x^{-4}),
\end{equation*}
where $f_B$ and $\beta_B$ means the density or $\beta$ of (B)
representation. We obtain (\ref{expansion}) by using the relation 
\begin{equation*}
f(x;\alpha,\beta)=\left(\cos\left(\frac{\pi}{2}\alpha\varrho\right)\right)^{\frac{1}{\alpha}}
f_B\left(\left(\cos\left(\frac{\pi}{2}\alpha\varrho\right)\right)^{\frac{1}{\alpha}}(x-\zeta);\alpha,\beta_B \right)
\end{equation*}
by noting that $\cos(\frac{\pi}{2}\alpha\varrho)\to 1$ and
$\beta_B=\beta+O(\Delta^2)$ as $\Delta \to 0$.
Then the problem is which of the two approximations (\ref{normal}) and (\ref{expansion})
is dominant for large $x$. We give an answer to this problem in the
following two theorems, which are based on Theorem 1 of Nagaev and Shkol'nik (1988). 
We can also see the precise statement on p.129 of Uchaikin and Zolotarev (2003) concerning this problem.

\begin{thm}
\label{thm:density-general-stable-distributions}
Let $\beta,\ -1<\beta<1$, be fixed and let $\beta^{\ast}=\beta\sgn(x-\zeta)$. Define 
\begin{eqnarray*}
F_1(x;\alpha,\beta) &=&f(x-\zeta;2), \\
F_2(x;\alpha,\beta) &=&(1+\beta^\ast)\Delta(x-\zeta)^{\Delta-3}, \\
g(x;\alpha,\beta) &=& F_1(x;\alpha,\beta)+F_2(x;\alpha,\beta). \label{eq:density-general-stable-distributions}
\end{eqnarray*}
For an arbitrarily small $\epsilon>0$ there
 exist $\Delta_0$ and $x_0$ such that for all $\Delta<\Delta_0$
 and $|x|>x_0$,
\begin{equation*}
\left|f(x;\alpha,\beta)/g(x;\alpha,\beta)-1\right|<\epsilon.
\end{equation*}
Furthermore, for an arbitrarily small constant $\delta>0$,
\begin{eqnarray*}
g(x;\alpha,\beta)=
\left\{
\begin{array}{ll}
F_1(x;\alpha,\beta)\left(1+o(\Delta^{\delta/2})\right) & \mbox{\rm if}\ |x-\zeta | \le (2-\delta)(\log1/\Delta)^{1/2}, \\ 
& \\
F_2(x;\alpha,\beta)\left(1+o(\Delta^{\delta/2})\right) & \mbox{\rm if}\ |x-\zeta | \ge (2+\delta)(\log1/\Delta)^{1/2}, 
\end{array}
\right.
\end{eqnarray*} 
uniformly in $|x|>x_0$.
\end{thm}

\begin{thm}
\label{thm:derivative-density}
Define 
\begin{eqnarray*}
F'_1(x;\alpha,\beta) &=&-\frac{x-\zeta}{2}f(x-\zeta;2), \\
F'_2(x;\alpha,\beta) &=&-3(1+\beta^\ast)\Delta(x-\zeta)^{\Delta-4}, \\
g'(x;\alpha,\beta) &=& F'_1(x;\alpha,\beta)+F'_2(x;\alpha,\beta). 
\end{eqnarray*}
Under the same conditions and notations of Theorem \ref{thm:density-general-stable-distributions}, for an arbitrarily small 
$\epsilon>0$ there
 exist $\Delta_0$ and $x_0$ such that for all $\Delta<\Delta_0$
 and $|x|>x_0$,
\begin{equation*}
\left|f'(x;\alpha,\beta)/g'(x;\alpha,\beta)-1\right|<\epsilon.
\end{equation*}
Furthermore, for an arbitrarily small constant $\delta>0$,
\begin{eqnarray*}
g'(x;\alpha,\beta)=
\left\{
\begin{array}{ll}
F'_1(x;\alpha,\beta)\left(1+o(\Delta^{\delta/2})\right) & \mbox{\rm if}\ |x-\zeta | \le (2-\delta)(\log1/\Delta)^{1/2}, \\ 
& \\
F'_2(x;\alpha,\beta)\left(1+o(\Delta^{\delta/2})\right) & \mbox{\rm if}\ |x-\zeta | \ge (2+\delta)(\log1/\Delta)^{1/2}, 
\end{array}
\right.
\end{eqnarray*} 
uniformly in $|x|>x_0$.
\end{thm}

We need the following 4 lemmas concerning behavior of
$A(\varphi;\alpha,\beta)$ in (\ref{integrand-density-A1}) for $\varphi \doteq 1$ to prove Theorem
\ref{thm:density-general-stable-distributions} and Theorem
\ref{thm:derivative-density}, which correspond to the lemmas of
Theorem 1 of Nagaev and Shkol'nik (1988).
In the lemmas and the proofs of Theorems
\ref{thm:density-general-stable-distributions} and \ref{thm:derivative-density}
$\epsilon$ and $\epsilon'$ denote arbitrarily small positive
constants because in the end of the proofs of Theorems
\ref{thm:density-general-stable-distributions} and \ref{thm:derivative-density}
we let $\epsilon$ and $\epsilon'$ $\downarrow 0$ separately. From now on
the Fisher information is evaluated at a fixed value of $\beta$,
$|\beta|<1$. For notational abbreviation we sometimes write 
$A(\varphi;\alpha,\beta)$ as $A(\varphi)$. We also write 
\begin{equation*}
\lambda=1-\varphi,\quad \varphi_\Delta=1-\Delta^{1/2-\epsilon}.
\end{equation*}

\begin{lem}
\label{lem:1-A}
 As $\Delta \to 0$ and for $0\le \lambda \le \Delta/\epsilon' $,
\begin{equation*}
\label{eq:lemma1-A}
A(1-\lambda)
=\frac{(\lambda/\Delta)^\frac{1}{1-\Delta}(1+\beta+\lambda/\Delta)}{(1+\beta+2\lambda/\Delta)^2}(1+o(\Delta\log(1/\Delta)))
\end{equation*}
uniformly in $\lambda$.
\end{lem}

\begin{lem}
\label{lem:2-A}
As $\Delta \to 0$
\begin{equation*}
A(\varphi_\Delta)=\frac{1}{4}+\frac{\pi^2}{16}\Delta^{1-2\epsilon}+o(\Delta),
 \quad
A'(\varphi_\Delta)=-\frac{\pi^2}{8}\Delta^{1/2-\epsilon}-\frac{(1+\beta)^2}{8}\Delta^{1/2+3\epsilon}+O(\Delta).
\end{equation*}
\end{lem}

\begin{lem}
\label{lem:3-A}
As $\Delta \to 0$ and $\lambda \to 0$ such that
 $\lambda\ge\Delta^{1/2-\epsilon}$, 
\begin{equation*}
A(1-\lambda)=\frac{1}{4}+\frac{\pi^2}{16}\lambda^2+o(\lambda^2), 
\quad 
A'(1-\lambda)=-\frac{\pi^2}{8}\lambda-\frac{\Delta^2(1-\lambda+\beta)^2}{8\lambda^3}+o(\lambda^2),
\end{equation*}
\begin{equation*}
A''(1-\lambda)=\frac{\pi^2}{8}+\frac{\pi^2}{8}\lambda^2-\frac{3}{8}\frac{\Delta^2(1-\lambda+\beta)^2}{\lambda^4}-
 \frac{3}{4}\frac{\Delta^2(1-\lambda+\beta)}{\lambda^3}\left(1-\frac{\Delta(1-\lambda+\beta)^2}{\lambda^2}  \right)+o(\lambda^2).
\end{equation*}
\end{lem}
Note that $\Delta/\epsilon' \le \Delta^{1/2-\epsilon} \le \lambda$ implies $
\Delta=o(\lambda^2)$ and we have to consider terms like $\Delta^2/\lambda^4$ or $\Delta^2/\lambda^3$
in Lemma \ref{lem:3-A}.
\begin{lem}
\label{lem:4-A}
If $\Delta$ is sufficiently small, $A(\varphi)$ is a monotonically
 decreasing function on $(-\varrho^{\ast},1)$.
\end{lem}
The proofs of the lemmas are given in Appendix.

In the following $\overline{c}>0$ is an appropriate positive constant.
\begin{pft}
We prove the assertion in the case $x-\zeta>0$ first.
As in the proof of Theorem 1 of Nagaev and Shkol'nik (1988), let $\tau,
 \eta>0$ be arbitrarily small numbers. We denote $z=(x-\zeta)^{\alpha/(\alpha-1)}$ and $\varphi_0=\varphi_\Delta-z^{-1/2+\tau}$.
Then we divide the integral
\begin{equation*}
H = \int^{1}_{-\varrho^{\ast}}   A(\varphi;\alpha,\beta)\exp(-zA(\varphi;\alpha,\beta)) d\varphi
\end{equation*}
into six subintegrals,
\begin{equation*}
H=\sum_{k=1}^{6}H_k,
\end{equation*}
where each $H_k$ corresponds to the integration of $H$ for the $k$-th interval
of $[-\varrho^{\ast}, 1-\eta)$, $[1-\eta, \varphi_0)$, $[\varphi_0, \varphi_\Delta)$, 
$[\varphi_\Delta, 1-\Delta/\epsilon')$, $[1-\Delta/\epsilon', 1-\Delta\epsilon')$
and $[1-\Delta\epsilon',1]$.
First we calculate $H_i's$ utilizing Lemmas \ref{lem:1-A}-\ref{lem:4-A}
 for fixed $\epsilon$ and $\epsilon'$. Later we see that $H_3$ and $H_6$ dominate the 
 others terms. We also let $\epsilon$ and $\epsilon'$ $\downarrow 0$ in the end.\vspace{2mm} \\ 
{\bf Calculation of} $H_1$: \\
From Lemma \ref{lem:2-A}, $A(\varphi_\Delta)>\frac{1}{4}$. By Lemma
 \ref{lem:4-A} and $1-\eta<\varphi_\Delta$, for sufficiently small
 $\Delta$ there exists a constant $\rho_1 \in (0,1)$ such that $\rho_1
 A(1-\eta)>\frac{1}{4}$. Then for a constant $\gamma>\frac{1}{4}$ we easily find 
\begin{eqnarray}
\label{eq:integration-H_1}
H_1 &=&
 \frac{1}{z(1-\rho_1)}\int^{1-\eta}_{-\varrho^\ast}z(1-\rho_1)A(\varphi)\exp\{-z\rho_1
 A(\varphi)-z(1-\rho_1)A(\varphi)\} d\varphi  \\
    &\le& \frac{1}{z(1-\rho_1)}\exp (-z \rho_1 A(1-\eta)) \nonumber \\
    &=& O(\exp(-\gamma z)/z). \nonumber
\end{eqnarray}   
\vspace{2mm} \\
{\bf Calculation of} $H_2$: \\
From Lemma \ref{lem:4-A}
\begin{eqnarray*}
H_2&=&\int^{\varphi_0}_{1-\eta} A(\varphi)\exp(-zA(\varphi)) d\varphi \\
   &\le&\eta A(1-\eta)\exp(-zA(\varphi_0)) \\
   &=& O(\exp(-zA(\varphi_0))).
\end{eqnarray*}
Then we can write 
\begin{equation*}
A(\varphi_0)=A(\varphi_\Delta)+A'(\varphi_\Delta)(\varphi_0-\varphi_\Delta)+\frac{1}{2}A''(\xi)(\varphi_0-\varphi_\Delta)^2,
\end{equation*}
where $\varphi_0 \le \xi \le \varphi_\Delta$. By Lemma \ref{lem:2-A} and
 Lemma \ref{lem:3-A} we can write
\begin{eqnarray*}
zA(\varphi_\Delta)&=&\frac{z}{4}+R_1(\Delta,z,\epsilon),  \\
zA'(\varphi_\Delta)(\varphi_0-\varphi_\Delta) &=& R_2(\Delta,z,\epsilon), \\
zA''(\xi)(\varphi_0-\varphi_\Delta)^2 &\le& \overline{c} z^{2\tau},
\end{eqnarray*}
where 
\begin{eqnarray*} 
R_1(\Delta,z,\epsilon) &=& O(z\Delta^{1-2\epsilon})>0,\\
R_2(\Delta,z,\epsilon) &=& O(\Delta^{1/2-\epsilon}z^{1/2+\tau})>0.
\end{eqnarray*}
Here we use functions $R_i(\Delta,z,\epsilon),\ i=1,2,$ for convenience in comparing $H_2$ with $H_3$. Then
\begin{equation}
\label{eq:integration-H_2}
H_2=O\left(\exp\left(-\frac{z}{4}-R_1(\Delta,z,\epsilon)-R_2(\Delta,z,\epsilon)-\overline{c}z^{2\tau}\right)\right).
\vspace{2mm}
\end{equation}
{\bf Calculation of} $H_3$: \\
We show the detailed calculation of $H_3$ because this is the dominant
 term in $H$ as $\Delta \to 0$ and $x \to \infty$.  
We have for $\xi \in [\varphi_0,\varphi_\Delta]$, 
\begin{equation*}
H_3=\exp(-zA(\varphi_\Delta))\int_{\varphi_0}^{\varphi_\Delta}A(\varphi)\exp\left(-zA'(\varphi_\Delta)(\varphi-\varphi_\Delta)
-\frac{z}{2} A''(\xi)(\varphi-\varphi_\Delta)^2\right) d\varphi.
\end{equation*}
Utilizing Lemmas \ref{lem:2-A} and  \ref{lem:3-A}, for some constants
 $\rho_1,\rho_2 \in (0,1)$ we obtain
\begin{eqnarray*}
A(\varphi) &=& \frac{1}{4}+\frac{\pi^2}{16}\lambda^2+o(\lambda^2), \\
-zA(\varphi_\Delta) &=&
-\frac{z}{4}-\frac{\pi^2}{16}\Delta^{1-2\epsilon}z+o(\Delta z), \\
-zA'(\varphi_\Delta)(\varphi-\varphi_\Delta) &=&
-\left(\frac{\pi^2}{8}\Delta^{1/2-\epsilon}+\frac{1}{8}(1+\beta)^2\Delta^{1/2+3\epsilon}\right)
\rho_1 z^{1/2+\tau}+O(\Delta z^{1/2+\tau}), \\
-\frac{z}{2}A''(\xi)(\varphi-\varphi_\Delta)^2
&=&
-\frac{z}{2}\frac{\pi^2}{8}(\varphi-\varphi_0)^2 \\
&&\hspace{-2cm}
-\left(\frac{\pi^2}{8}\lambda^2-\frac{3}{8}\frac{\Delta^2(1-\lambda+\beta)^2}{\lambda^4}
-\frac{3}{4}\frac{\Delta^2(1-\beta+\lambda)}{\lambda^3}\left(1-\frac{\Delta(1-\lambda+\beta)^2}{\lambda^2} \right)+o(\lambda^2)
\right)\frac{\rho_2}{2} z^{2\tau}. 
\end{eqnarray*}
Summarizing the main terms, we obtain
\begin{eqnarray*}
A(\varphi) &=& \frac{1}{4}+O(\lambda^2), \\
-zA(\varphi_\Delta) &=&
-\frac{z}{4}-R_1(\Delta,z,\epsilon), \\
-zA'(\varphi_\Delta)(\varphi-\varphi_\Delta) &=&
-R_3(\Delta,z,\epsilon), \\
-\frac{z}{2}A''(\xi)(\varphi-\varphi_\Delta)^2
&=&
-\frac{z}{2}\frac{\pi^2}{8}(\varphi-\varphi_0)^2+O(\lambda^2z^{2\tau})+O(\Delta^2/\lambda^4z^{2\tau})\\ 
&=&
-\frac{z}{2}\frac{\pi^2}{8}(\varphi-\varphi_0)^2+O(z^{-1+4\tau})+O(\Delta^{4\epsilon}z^{2\tau}),
\end{eqnarray*}
where 
\begin{equation*}
0<R_3(\Delta,z,\epsilon)=O(\Delta^{1/2-\epsilon}z^{1/2+\tau}) \le R_2(\Delta,z,\epsilon).
\end{equation*}
Note $\varphi_0 \le \varphi \le \varphi_\Delta  \Leftrightarrow \ \Delta^{1/2-\epsilon} \le \lambda \le \Delta^{1/2-\epsilon}+z^{-1/2+\tau}$.
Then we obtain
\begin{equation*}
H_3=R(\lambda,\Delta,z,\epsilon)\frac{1}{4}\exp\left(-\frac{z}{4}\right)
\int_{\varphi_0}^{\varphi_\Delta}\exp\left(-\frac{z}{2}\frac{\pi^2}{8}
(\varphi-\varphi_\Delta)^2 \right) d\varphi, 
\end{equation*}
where 
\begin{equation*}
R(\lambda,\Delta,z,\epsilon) = \left(1+O(\lambda^2) \right)
\exp\left(-R_1(\Delta,z,\epsilon)-R_3(\Delta,z,\epsilon)+O(z^{-1+4\tau})+O(\Delta^{4\epsilon}z^{2\tau})\right).
\end{equation*}
Here the inequality
\begin{equation*}
\int_x^{\infty}e^{-y^2/2}dy \le \frac{1}{x}e^{-x^2/2}, \quad x>0
\end{equation*}
gives
\begin{eqnarray*}
\int_{\varphi_0}^{\varphi_\Delta}\exp\left(-\frac{z}{2}\frac{\pi^2}{8}
(\varphi-\varphi_0)^2\right) d\varphi &=&
\frac{1}{\sqrt{z}}\int_{0}^{\infty}\exp\left(-\frac{z}{2}\frac{\pi^2}{8}
\varphi^2 \right) d\varphi
-
\frac{1}{\sqrt{z}}\int_{z^{\tau}}^{\infty}\exp\left(-\frac{z}{2}\frac{\pi^2}{8}
\varphi^2\right) d\varphi \\
&=&
\frac{2}{\sqrt{\pi z}}+O\left(\frac{1}{z^{1/2+\tau}}\exp\left(\frac{-z^{2\tau}
}{2}\right)\right).
\end{eqnarray*}
Then we obtain
\begin{equation}
\label{eq:integration-H_3-1}
H_3=\exp\left(-\frac{z}{4}-R_1(\Delta,z,\epsilon)-R_3(\Delta,z,\epsilon)\right)\frac{1}{2\sqrt{\pi
 z}}\exp\left(O(z^{-1+4\tau})+O(\Delta^{4\epsilon}z^{2\tau}) \right)(1+o(1)).
\end{equation}
Written in this form, $H_3$ can be easily compared with $H_2$, i.e., the
 formula (\ref{eq:integration-H_2}). Further
 we will see later that $H_3$ will be dominant if $z \le (4-\delta)\log1/\Delta $ for an arbitrarily small constant
 $\delta>0$. Then we only have to consider the convergence of $H_3$ for $z \le O(\log1/\Delta)$.
 For $z \le O(\log1/\Delta)$,
 $R_1(\Delta,z,\epsilon) \downarrow 0$,
 $R_2(\Delta,z,\epsilon) \downarrow 0$ and $R_3(\Delta,z,\epsilon)
 \downarrow 0$ as
 $\Delta \to 0$.
Finally we obtain
\begin{equation}
\label{eq:integration-H_3}
H_3=\frac{1}{2\sqrt{\pi z}}\exp\left(-\frac{z}{4}\right)(1+o(1)) \quad \mbox{if}\quad z \le O(\log1/\Delta). 
\vspace{2mm}
\end{equation}
{\bf Calculation of} $H_4$:\\
From Lemma \ref{lem:4-A}
\begin{eqnarray*}
H_4 
&\le& \int^{1-\Delta/\epsilon'}_{\varphi_\Delta}
A(\varphi_\Delta)\exp(-zA(1-\Delta/\epsilon')) d\varphi \\
&\le& \Delta^{1/2-\epsilon}A(\varphi_\Delta)\exp(-zA(1-\Delta/\epsilon')).
\end{eqnarray*}
Here $\lambda/\Delta \le 1/\epsilon'$. Then by Lemma \ref{lem:1-A} and
 \ref{lem:4-A} as $\Delta \to 0$
\begin{equation*}
A(1-\Delta/\epsilon')=\frac{{\epsilon'}^{\frac{\Delta}{\Delta-1}}\{\epsilon'(1+\beta)+1\}}{\{\epsilon'(1+\beta)+2\}^2}
=\frac{1}{4}(1-w(\epsilon')),
\end{equation*}
where $w(\cdot)$ is a nonnegative function such as
\begin{equation*}
\lim_{t \to 0}w(t)=0.
\end{equation*}
Finally we obtain
\begin{equation}
\label{eq:integration-H_4}
H_4 \le O\left(\Delta^{1/2-\epsilon}\exp\left(-\frac{z}{4}(1-w(\epsilon'))\right)\right).
\end{equation}
For the comparison of $H_4$ with $H_3$ and $H_6$, which is needed in
later argument, we state some properties of $H_4$.
For a given small number $\delta>0$, there exists $\epsilon'$ such that
$\delta>w(\epsilon')$.
Then for $z \le 4\log1/\Delta$,
\begin{equation*}
\Delta^{\delta}=\exp\left(-\frac{\delta}{4}\times4\log1/\Delta\right) \le \exp\left(-\frac{\delta}{4}z\right).
\end{equation*}
Consequently
\begin{equation}
\label{eq:integration-H_4-1}
H_4 \le O\left(\Delta^{1/2-\delta-\epsilon} \exp\left(-\frac{z}{4}(1+\delta-w(\epsilon'))\right)\right).
\end{equation}
In this representation $H_4$ can be easily compared with $H_3$.
For $z\ge 4\log1/\Delta$
\begin{equation}
\label{eq:integration-H_4-2}
H_4 \le O(\Delta \exp(-\overline{c}z)).
\end{equation}
In this representation $H_4$ can be easily compared with $H_6$.
\vspace{2mm}\\
{\bf Calculation of} $H_5$: \\
From Lemma \ref{lem:1-A} for $\epsilon' \le \lambda/\Delta \le 1/\epsilon'$,
$A(\varphi)$ is bounded. Then
\begin{eqnarray}
\label{eq:integration-H_5}
H_5 &=& \int^{1-\Delta\epsilon'}_{1-\Delta/\epsilon'}A(\varphi)
\exp(-zA(\varphi)) d\varphi \\
&\le& (\Delta/\epsilon'-\Delta\epsilon')\sup_{\varphi \in (1-\Delta/\epsilon',1-\Delta\epsilon')}
 A(\varphi)\exp(-zA(\varphi)) \nonumber \\
&\le&
O(\Delta\exp(-\overline{c}z)). \nonumber
\vspace{2mm}
\end{eqnarray}   
{\bf Calculation of} $H_6$:\\
From Lemma \ref{lem:1-A} for $\lambda/\Delta \le \epsilon'$
\begin{eqnarray*}
A(1-\lambda;\alpha,\beta) 
&=& 
\frac{(\lambda/\Delta)^{1/(1-\Delta)}}{1+\beta}\left(1+O(\epsilon')+O(\Delta\log(1/\Delta))\right)\\
&=&
\frac{(\lambda/\Delta)^{1/(1-\Delta)}}{1+\beta}(1+R_4(\Delta,\epsilon')),
\end{eqnarray*}
where $R_4(\Delta,\epsilon')=O(\epsilon')+O(\Delta\log(1/\Delta))$. Substituting the above results
into $H_6$ and replacing 
$\varphi=1-\lambda$ by $\lambda$, we obtain
\begin{equation*}
H_6=\left(\frac{1+R_4}{1+\beta}\right)
\int_0^{\Delta\epsilon'}\left(\lambda/\Delta\right)^{1/(1-\Delta)}
\exp\left(-z\left(\frac{1+R_4}{1+\beta}\right)\left(\lambda/\Delta\right)^{1/(1-\Delta)}
\right) d\lambda.
\end{equation*}
Furthermore, replacing 
$$
\lambda=\Delta \left(\frac{1+\beta}{1+R_4}\frac{x}{z}\right)^{1-\Delta}
$$
by $x$ and defining  
$$
g(z)={\epsilon'}^{1/(1-\Delta)}\left(\frac{1+R_4}{1+\beta}\right)z,
$$
$H_6$ is written as
\begin{equation*}
H_6 =
\frac{\Delta(1+\beta)^{1-\Delta}}{z^{2-\Delta}}
\int_0^{g(z)}
x\exp(-x)x^{-\Delta}dx \times \{(1+R_4)^{(\Delta-1)}(1-\Delta)\} . 
\end{equation*}
Then the integration of the right side of the equation is evaluated as
\begin{eqnarray*}
\int_0^{g(z)}
x \exp(-x)dx+O(\Delta) 
&=&
\left[ -xe^{-x}-e^{-x}
\right]_0^{g(z)}+O(\Delta) \\
&=&
1+O\left( g e^{-g}\right)
+O(\Delta).
\end{eqnarray*}
 Here for arbitrarily small constant $\epsilon'$ and any $\Delta \to 0$, $g(z) \to \infty$ as $z \to \infty$. Therefore
\begin{equation}
\label{eq:integration-H_6}
H_6 = \Delta(1+\beta) z^{\Delta-2}(1+o(1)).
\end{equation} 
For fixed $\epsilon$ and $\epsilon'$ the formulas (\ref{eq:integration-H_1}), (\ref{eq:integration-H_2}), (\ref{eq:integration-H_3-1}), (\ref{eq:integration-H_4}), 
(\ref{eq:integration-H_5}) and (\ref{eq:integration-H_6}) are proved.

Comparing (\ref{eq:integration-H_1}), (\ref{eq:integration-H_2}),
(\ref{eq:integration-H_3-1}), (\ref{eq:integration-H_4}), 
(\ref{eq:integration-H_5}) and (\ref{eq:integration-H_6}), the order of
 $H_i, i=1,\ldots,6$,
 is summarized as follows: 
\begin{equation*}
\mbox{max}(H_1,H_2,H_4,H_5)=o(\mbox{max}(H_3,H_6)).
\end{equation*}
The relative dominance of $H_3$ (\ref{eq:integration-H_3}) and $H_6$ (\ref{eq:integration-H_6}) depends on the value of $z$.
For an arbitrarily small $\delta>0$,
\begin{eqnarray*}
H_6=o(H_3)  \quad \mbox{if}\ z \le  (4-\delta)\log 1/\Delta, \\
H_3=o(H_6)  \quad \mbox{if}\ z \geq (4+\delta)\log 1/\Delta.
\end{eqnarray*}
These are easily confirmed if we substitute $(4\pm\delta)\log 1/\Delta$
into $z$.

Now for proving the uniformity in $\Delta$ and $z$  we let $\epsilon,
 \epsilon' \downarrow 0$ depending on $\Delta$ and $z$. We first let $\epsilon \downarrow 0$. We can take
 $\epsilon=(\log1/\Delta)^{\delta-1}$ for a constant $\delta \in (0,1)$. Concerning Lemmas 
 \ref{lem:1-A}-\ref{lem:4-A}, we must consider the relations of
 $\Delta^{2\epsilon} \ge \Delta/\lambda^2$ in Lemma \ref{lem:3-A}. In
 the other lemmas the results do not change. Substituting
 $\epsilon=(\log1/\Delta)^{\delta-1}$ into $\Delta^{2\epsilon}$ and
 taking logarithm of $\Delta^{2\epsilon}$ 
\begin{equation*}
\log \Delta^{2\epsilon}=-2 (\log1/\Delta)^{\delta} \to -  \infty.
\end{equation*}
Hence $\Delta=o(\lambda^2)$ and the conclusion of Lemma \ref{lem:3-A} is
 not changed. Concerning the calculations of $H_i$, the order of dominance is not
 affected by the limiting operation $\epsilon \downarrow 0$. For the comparison of
 $H_3$ and $H_6$ we may substitute $\epsilon=(\log1/\Delta)^{\delta-1}$
 and $z=(4\pm\delta)\log 1/\Delta$ into (\ref{eq:integration-H_3-1}) and 
(\ref{eq:integration-H_6}). The other comparisons are easy. Then we have
 only to confirm that $H_3$ satisfies (\ref{eq:integration-H_3})
 because the calculation of $H_6$ is not affected by $\epsilon$. If we
 substitute $\epsilon=(\log1/\Delta)^{\delta-1}$ into
 (\ref{eq:integration-H_3-1}), we obtain (\ref{eq:integration-H_3}) since  
\begin{equation*}
z\Delta^{1-2\epsilon}\to 0,\quad \Delta^{1/2-\epsilon}z^{1/2+\tau}\to 0,\quad
 \Delta^{4\epsilon}z^{2\tau}\to 0,\quad \mbox{for}\ z \le 4\log1/\Delta
\end{equation*}
imply
\begin{equation*}
R_1(\Delta,z,\epsilon) \downarrow 0, \quad
R_2(\Delta,z,\epsilon) \downarrow 0, \quad
R_3(\Delta,z,\epsilon) \downarrow 0, \quad
O(\Delta^{4\epsilon}z^{2\tau}) \to 0, \quad \mbox{for}\ z \le 4\log1/\Delta.
\end{equation*}

Next we consider $\epsilon' \downarrow 0$ depending on $\Delta$. We may
 consider $\epsilon'=(\log1/\Delta)^{\delta'-1}$ for an arbitrarily small
 $\delta' \in (0,1)$. Concerning Lemma \ref{lem:1-A} $\epsilon'$ must satisfy $\Delta/\epsilon' \le
 \Delta^{1/2-\epsilon}$. This is easily confirmed. $\epsilon'$ does not
 affect the other lemmas. Concerning the calculations $H_i's$, we must
 consider $H_4$, $H_5$ and $H_6$. The other $H_i's$ are not affected by
 $\epsilon'$. From (\ref{eq:integration-H_4}) and the definition of
 $w(\epsilon')$, the calculation of $H_4$ concerning the order, i.e., formulas (\ref{eq:integration-H_4-1})
 and (\ref{eq:integration-H_4-2}) are not changed. We investigate $H_5$
 in detail. If $[1-\Delta/\epsilon,1-\Delta\epsilon')$, then $\epsilon'
 \le \lambda/\Delta \le 1/\epsilon'$ and $A(\varphi)$ of Lemma
 \ref{lem:1-A} satisfies 
\begin{equation*}
\frac{1+o(1)}{1+\beta} {\epsilon'}^{1/(1-\Delta)} \le A(\varphi).
\end{equation*}
Hence
\begin{eqnarray*}
H_5 &\le& (\Delta/\epsilon'-\Delta\epsilon') \times \sup_{\varphi \in
 (1-\Delta/\epsilon,1-\Delta\epsilon)} 
 A(\varphi)\exp(-zA(\varphi))  \\
    &\le& O\left( (\Delta/\epsilon'-\Delta\epsilon') \times
      {\epsilon'}^{\frac{1}{1-\Delta}} \exp
      \left(-\overline{c}z{\epsilon'}^{\frac{1}{1-\Delta}}\right)\right) \nonumber \\
     &\le& O\left(\Delta \exp\left(-\overline{c}z{\epsilon'}^{1/(1-\Delta)}\right)\right). \nonumber 
\end{eqnarray*}
Then for $z \le (4-4\delta'')\log1/\Delta$ and for an arbitrarily small
 positive constant
 $\delta''$, 
\begin{equation*}
H_5 \le O(\Delta) \le O\left(\Delta^{\delta''} \exp\left(-\frac{z}{4}\right)\right) = o(H_3).
\end{equation*}
For $z \ge (4-4\delta'')\log 1/\Delta$
\begin{equation*}
H_5=O(\Delta \exp(-\overline{c}z^{\delta'})) =o(H_6).
\end{equation*} 
Hence the order and dominant terms of $H_i$'s are not changed. For the
 calculation of $H_6$, $g(z)$ must go to $\infty$ as $z \to \infty$. Considering $z \ge
 (4+\delta)\log1/\Delta$, we can easily confirm this condition. Thus we
 can let $\epsilon' \downarrow 0$ depending on $\Delta$. 
Therefore convergence of $H_i,i=1,\ldots,6$, only depends on
$\Delta$ and $z$ and is faster as $\Delta \to 0$ and $z \to \infty$.

Finally, substituting $z=(x-\zeta)^{\alpha/(\alpha-1)}$ into $H_3$ and $H_6$ and 
multiplying these by $\alpha/(2(\alpha-1))(x-\zeta)^{1/(\alpha-1)}$, we
obtain the result for $x-\zeta>0$. For $x-\zeta<0$ we utilize the relations
$f(x;\alpha,\beta)=f(-x;\alpha,-\beta)$ of general stable
distributions. 
 \hfill $\Box$ \\
\end{pft}
\begin{pft}{}
From (\ref{eq:dense})
the derivative of the density is calculated as follows.
\begin{eqnarray}
&& f'(x;\alpha,\beta)=\label{eq:derivatives-of-density}
\frac{1}{\alpha-1}(x-\zeta)^{-1} f(x;\alpha,\beta)-\frac{\alpha^2}{2(\alpha-1)^2}(x-\zeta)^{2/(\alpha-1)}\times
 H', 
\end{eqnarray}
where
\begin{equation*}
H'=\int_{-\varrho^\ast}^{1}
 A^2(\varphi;\alpha,\beta)\exp\left(-(x-\zeta)^{\alpha/(\alpha-1)}A(\varphi;\alpha,\beta)\right) d\varphi. 
\end{equation*}
Since behavior of $f(x;\alpha,\beta)$ was obtained in Theorem
 \ref{thm:density-general-stable-distributions}, we evaluate the 
 integration of the second term of (\ref{eq:derivatives-of-density}).
As in Theorem \ref{thm:density-general-stable-distributions} we divide
 integral $H'$
into subintegrals,
\begin{equation*}
H'=\sum_{k=1}^{6}H'_k,
\end{equation*}
where each $H'_k$ corresponds to the $k$-th interval
of $[0, 1-\mu)$, $[1-\mu, \varphi_0)$, $[\varphi_0, \varphi_\Delta)$, 
$[\varphi_\Delta, 1-\Delta/\epsilon')$, $[1-\Delta/\epsilon', 1-\Delta\epsilon')$
and $[1-\Delta\epsilon',1]$. From Lemmas \ref{lem:1-A}, \ref{lem:2-A} and
 \ref{lem:3-A},  
for $\lambda \le \Delta/\epsilon' $
\begin{equation}
\label{eq:lemma1-A^2}
A^2(1-\lambda)=\frac{(\lambda/\Delta)^\frac{2}{1-\Delta}(1+\beta+\lambda/\Delta)^2}{(1+\beta+2\lambda/\Delta)^4}(1+O(\Delta\log(1/\Delta)))
, \\
\end{equation}
\begin{equation*}
A^2(\varphi_\Delta)=\frac{1}{16}+\frac{\pi^2}{32}\Delta^{1-2\epsilon}+o(\Delta),
\end{equation*}
and for $\varphi \le \varphi_\Delta \Leftrightarrow
 \Delta^{1/2-\epsilon} \le \lambda$
\begin{equation*}
A^2(1-\lambda)=\frac{1}{16}+\frac{\pi^2}{32}\lambda^2+o(\lambda^2).
\end{equation*}
Considering $A^2(\varphi,\Delta)$ term of the integrand, we easily find the
 following results on $H'_i$ similar to $H_i, i=1,\ldots,5,$ in the proof of Theorem \ref{thm:density-general-stable-distributions}.
\begin{eqnarray}
\label{eq:integration-H'_1}
H'_1 &\le& O(\exp(-\gamma z)/z^2), \quad \mbox{for} \quad \gamma >\frac{1}{4},\\
\label{eq:integration-H'_2}
H'_2 &\le& O\left(\exp \left(-\frac{z}{4}+O(z\Delta^{1-2\epsilon})
+O(z^{1/2+\tau}\Delta^{1/2-\epsilon})-\overline{c}z^{2\tau}\right)\right), \\
\label{eq:integration-H'_3}
H'_3 &=&\frac{1}{8\sqrt{z\pi}}\exp\left(-\frac{z}{4}\right)(1+o(1)),\\
\label{eq:integration-H'_4}
H'_4 &\le& O\left(\Delta^{1/2-\epsilon}\exp\left(-\frac{z}{4}(1-w(\epsilon'))\right)\right), \\
\label{eq:integration-H'_5}
H'_5 &\le& O(\Delta\exp(-\overline{c}z)).
\end{eqnarray}
Later we see $H'_3$ and $H'_6$ dominate the 
 others terms as $H_3$ and $H_6$ dominate the other $H_i$'s in Theorem \ref{thm:density-general-stable-distributions}.
Since the calculation $H'_6$ is somewhat different, we treat $H'_6$ separately.
From (\ref{eq:lemma1-A^2}), for $\lambda/\Delta \le \epsilon'$
\begin{eqnarray*}
A^2(1-\lambda;\alpha,\beta) 
&=& 
\frac{(\lambda/\Delta)^\frac{2}{1-\Delta}}{(1+\beta)^2}\left(1+R_4(\Delta,\epsilon')\right)^2,
\end{eqnarray*}
where $R_4(\Delta,\epsilon')$ is defined in the proof of Theorem \ref{thm:density-general-stable-distributions}. Substituting the above
 result into $H'_6$ and replacing 
$\varphi=1-\lambda$ by $\lambda$, we obtain
\begin{equation*}
H'_6=\left(\frac{1+R_4}{1+\beta}\right)^2
\int_0^{\Delta\epsilon'}\left(\lambda/\Delta\right)^{2/(1-\Delta)}
\exp\left(-z\left(\frac{1+R_4}{1+\beta}\right)\left(\lambda/\Delta\right)^{1/(1-\Delta)}
\right) d\lambda.
\end{equation*}
Furthermore, replacing 
$$
\lambda=\Delta \left(\frac{1+\beta}{1+R_4}\frac{x}{z}\right)^{1-\Delta}
$$
by $x$ and defining 
$$
g(z)={\epsilon'}^{1/(1-\Delta)}\left(\frac{1+R_4}{1+\beta}\right)z,
$$
$H'_6$ is written as 
\begin{equation*}
H'_6 =
\frac{\Delta(1+\beta)^{1-\Delta}}{z^{3-\Delta}}
\int_0^{g(z)}
x^2\exp(-x)x^{-\Delta}dx \times \{(1+R_4)^{(\Delta-1)}(1-\Delta)\} . 
\end{equation*}
Then the integration of the right side of the equation is evaluated as
\begin{eqnarray*}
\int_0^{g(z)}
x^2 \exp(-x)dx+O(\Delta) 
&=&
\left[ -x^2e^{-x}-2xe^{-x}-2e^{-x}
\right]_0^{g(z)}+O(\Delta) \\
&=&
2+O\left( g^2 e^{-g}\right) +O(\Delta).
\end{eqnarray*}
 Here for arbitrarily small constant $\epsilon'$ and any $\Delta \to 0$, $g(z) \to \infty$ as $z \to \infty$. Therefore
\begin{equation}
\label{eq:integration-H'_6}
H'_6 = 2(1+\beta)\Delta z^{\Delta-3}(1+o(1)).
\end{equation}
Hence for fixed $\epsilon$ and $\epsilon'$ formulas (\ref{eq:integration-H'_1})-(\ref{eq:integration-H'_6}) are proved. 

From (\ref{eq:integration-H'_1})-(\ref{eq:integration-H'_6}) 
$$
\mbox{max}(H'_1,H'_2,H'_4,H'_5)=o (\mbox{max}(H'_3,H'_6)).
$$
The relative dominance $H_3'$ and $H_6'$ depend on the value of $z$.
For an arbitrarily small $\delta>0$    
\begin{eqnarray*}
H'_6=o(H'_3) \quad  \mbox{if}\ z \le  (4-\delta)\log1/\Delta, \\
H'_3=o(H'_6) \quad  \mbox{if}\ z \geq (4+\delta)\log1/\Delta.
\end{eqnarray*}
These are easily confirmed if we substitute $(4\pm\delta)\log 1/\Delta$
into $z$.

Substituting $z=(x-\zeta)^{\alpha/(\alpha-1)}$ into $H'_3$ and $H'_6$ and
multiplying each of them by $\alpha^2/(2(\alpha-1)^2)(x-\zeta)^{2/(\alpha-1)}$ 
we obtain the second integral of the equation
 (\ref{eq:derivatives-of-density}). 
Combining this with Theorem \ref{thm:density-general-stable-distributions}, we obtain the desired result. For $x-\zeta<0$ we use the relation
$f(x;\alpha,\beta)=f(-x;\alpha,-\beta)$.
For proving uniformity in
 $\Delta$ and $z$ we need to let $\epsilon, \epsilon'
 \downarrow 0$.
We omit the derivation because it is quite similar to that in the proof of Theorem
 \ref{thm:density-general-stable-distributions}. \\
 \hfill $\Box$ 
\end{pft}
From Theorem \ref{thm:density-general-stable-distributions}, Theorem
\ref{thm:derivative-density} and formulas (\ref{eq:f-mu}),
(\ref{eq:f-sigma}), we obtain  
behavior of the score functions as $\alpha \uparrow 2$. These score
functions in the corollaries below are needed
for obtaining 
the Fisher information matrix of the general stable distributions. 
\begin{cor}
\label{cor:score-location}
Under the same conditions and notations of Theorem
 \ref{thm:density-general-stable-distributions} and \ref{thm:derivative-density}, define
\begin{equation*}
g_\mu(x;\alpha,\beta)=-g'(x;\alpha,\beta)/g(x;\alpha,\beta).
\end{equation*}
Then for an arbitrarily small 
$\epsilon>0$, there
 exist $\Delta_0$ and $x_0$ such that for all $\Delta<\Delta_0$
 and $|x|>x_0$,
\begin{equation*}
|f_\mu(x;\alpha,\beta)/g_\mu(x;\alpha,\beta)-1|<\epsilon.
\end{equation*}
Furthermore define
\begin{eqnarray*}
F_\mu^1(x;\alpha,\beta)&=&
\frac{x-\zeta}{2}, \\
F_\mu^2(x;\alpha,\beta)&=&\frac{3}{x-\zeta}.
\end{eqnarray*}
Then 
 \begin{eqnarray*}
\frac{g_\mu(x;\alpha,\beta)}
{g(x;\alpha,\beta)}
=
\left\{
\begin{array}{ll}
F_\mu^1(x;\alpha,\beta)(1+o(\Delta^{\delta/2})) & \mbox{\rm if}\ |x-\zeta | \le (2-\delta)(\log1/\Delta)^{1/2} \\ 
& \\
F_\mu^2(x;\alpha,\beta)(1+o(\Delta^{\delta/2})) & \mbox{\rm if}\ |x-\zeta | \ge (2+\delta)(\log1/\Delta)^{1/2},
\end{array}
\right.
\end{eqnarray*} 
uniformly in $|x|>x_0$.
For the remaining interval $(2-\delta)(\log1/\Delta)^{1/2}\le x \le (2+\delta)(\log1/\Delta)^{1/2}$,
\begin{equation*}
\left| \frac{g_\mu(x;\alpha,\beta)}
{g(x;\alpha,\beta)}  \right| \le \overline{c} x
\end{equation*}
uniformly in $|x|>x_0$.
\end{cor}

\begin{cor}
\label{cor:score-scale}
Under the same conditions and notations of Theorems
 \ref{thm:density-general-stable-distributions} and \ref{thm:derivative-density}, define
\begin{equation*}
g_\sigma(x;\alpha,\beta)=-g(x;\alpha,\beta)-(x-\zeta)g'(x;\alpha,\beta).
\end{equation*}
Then for an arbitrarily small 
$\delta>0$, there
 exist $\Delta_0$ and $x_0$ such that for all $\Delta<\Delta_0$
 and $|x|>x_0$
\begin{equation*}
|f_\sigma(x;\alpha,\beta)/g_\sigma(x;\alpha,\beta)-1|<\epsilon.
\end{equation*}
Furthermore define
\begin{eqnarray*}
F_\sigma^1(x;\alpha,\beta)&=&
\frac{(x-\zeta)^2}{2}, \\
F_\sigma^2(x;\alpha,\beta)&=&2.
\end{eqnarray*}
Then 
 \begin{eqnarray*}
\frac{g_\sigma(x;\alpha,\beta)}
{g(x;\alpha,\beta)}
=
\left\{
\begin{array}{ll}
F_\sigma^1(x;\alpha,\beta)(1+o(\Delta^{1/2})) & \mbox{\rm if}\ |x-\zeta | \le (2-\delta)(\log1/\Delta)^{1/2} \\ 
& \\
F_\sigma^2(x;\alpha,\beta)(1+o(\Delta^{1/2})) & \mbox{\rm if}\ |x-\zeta | \ge (2+\delta)(\log1/\Delta)^{1/2}. 
\end{array}
\right.
\end{eqnarray*} 
For the remaining interval $(2-\delta)(\log1/\Delta)^{1/2}\le x \le (2+\delta)(\log1/\Delta)^{1/2}$,
\begin{equation*}
\left| \frac{g_\sigma(x;\alpha,\beta)}
{g(x;\alpha,\beta)}  \right| \le \overline{c} x^2
\end{equation*}
uniformly in $|x|>x_0$.
\end{cor}  

\section{Derivatives of density with respect to the parameters $\alpha$ and $\beta$}
In this section we obtain the derivatives of the density with respect to
the parameters $\alpha$ and $\beta$ by
analyzing the inversion formula. The density of general stable distributions
 can be written as
\begin{equation}
\label{eq:density-inversion-formula-half}
f(x;\alpha,\beta)=
\frac{1}{\pi}\mbox{Re}\int^\infty_0 e^{-itx}\Phi(t;\alpha,\beta)dt,
\end{equation}
where $\Phi(t;\alpha,\beta)$ is defined in
(\ref{eq:characteristic-function-standard}).  
Utilizing (\ref{eq:density-inversion-formula-half}), we derive the following
two lemmas concerning derivatives with respect to $\alpha$ and $\beta$, which are
extensions of Lemma 5 of Nagaev  and  Shkol'nik (1988) to non-symmetric
case.  

\begin{lem}
\label{lem:f_alpha(x;alpha,beta)}
As $\Delta=2-\alpha \to 0$, there exists $x_0$ and for all $|x| \ge x_0$,
\begin{equation*}
f_\alpha(x;\alpha,\beta)=
-\frac{1}{|y|^{1+\alpha}}\left\{
1+\beta^\ast+\Delta (M_1+M_2\log |y|)+\frac{M_3}{|y|}+\frac{M_4+M_5\log |y|}{|y|^\alpha}
\right\}, 
\end{equation*}
where $y=x-\zeta$ and $M_1,\cdots,M_5$ are some constants.
\end{lem}

\begin{lem}
\label{lem:f_beta(x;alpha,beta)}
Under the same conditions and notations of Lemma \ref{lem:f_alpha(x;alpha,beta)},
\begin{equation*}
f_\beta(x;\alpha,\beta)=
\frac{\Delta\sgn y}{|y|^{1+\alpha}}\left(1+\Delta M_6+\frac{M_7}{|y|}+\frac{M_8}{|y|^{\alpha}}\right),
\end{equation*}
where $M_6$, $M_7$ and $M_8$ are some constants.
Furthermore
\begin{equation*}
|f_\beta(x;\alpha,\beta)|=O(\Delta),
\end{equation*}
uniformly $x \in \mbox{\bf R}$.
\end{lem}

\begin{pfl}
First we consider the
case of $y=x-\zeta>0$. 
Differentiating $f$ by $\alpha$, we get
\begin{eqnarray}
\label{eq:differential-alpha-density-inversion-formula-half}
f_\alpha(x;\alpha,\beta)&=& \frac{1}{\pi} \mbox{Re} \int^\infty_0
e^{-ity}\exp\left(-t^\alpha\left(1-i\beta\tan\left(\frac{\pi\alpha}{2}\right)\right)\right) \\
&&\times \left\{-t^\alpha\log
 t\left(1-i\beta\tan\left(\frac{\pi\alpha}{2}\right)\right)+(t^\alpha-t)\frac{i\beta\pi}{2}/\left(\cos\left(\frac{\pi\alpha}{2}\right)\right)^2\right\} dt. \nonumber
\end{eqnarray}
Transforming $f_\alpha(x;\alpha,\beta)$ by $t \to s/y$, we obtain
\begin{eqnarray*}
f_\alpha(x;\alpha,\beta)&=&
\frac{1}{\pi y}\mbox{Re}\int^\infty_0 e^{-is}\exp\left(-s^\alpha
						  y^{-\alpha}\left(1-i\beta\tan\left(\frac{\pi\alpha}{2}\right)\right)^2\right) \nonumber \\
&& \times\left\{-s^{\alpha}y^{-\alpha}(\log s+\log
 y)\left(1-i\beta\tan\left(\frac{\pi\alpha}{2}\right)\right)+(s^\alpha
 y^{-\alpha}-sy^{-1})\frac{i\beta\pi}{2}/\left(\cos\left(\frac{\pi\alpha}{2}\right)\right)^2\right\} ds \nonumber
\end{eqnarray*}
Further transforming $f_\alpha(x;\alpha,\beta)$ by $s \to te^{i\varsigma}, \varsigma=-\frac{\pi}{2\alpha}$
as in Lemma 5 of Nagaev  and  Shkol'nik (1988), we find
\begin{eqnarray*}
f_\alpha(x;\alpha,\beta)
&=&
\frac{1}{\pi y}\mbox{Re}\ e^{i\varsigma}
\int^{\infty}_{0}\exp\left(-ite^{i\varsigma}+it^{\alpha}y^{-\alpha}\left(1-i\beta\tan\left(\frac{\pi\alpha}{2}\right)\right)\right)
\nonumber \\
&&\times 
\left\{ it^{\alpha}y^{-\alpha}(\log t+i\varsigma+\log
 y)\left(1-i\beta\tan\left(\frac{\pi\alpha}{2}\right)\right)
-(it^{\alpha}y^{-\alpha}+ty^{-1}e^{i\varsigma})\frac{i\beta\pi}{2}/\left(\cos\frac{\pi\alpha}{2}\right)^2
\right\} dt \nonumber \\
&=& K_1+K_2, \nonumber
\end{eqnarray*}
where $K_1$ and $K_2$ are obtained in formulas (\ref{eq:K1}) and
 (\ref{eq:K2}) of appendix \ref{ap:proof-of-lemma3.1}.
 Then the lemma is proved for $y=x-\zeta>0$. 
For $y=x-\zeta<0$, we use the relation
$f_\alpha(x;\alpha,\beta)=f_\alpha(-x;\alpha,-\beta)$.\\
 \hfill $\Box$ \\
\end{pfl}
\begin{pfl}
First we assume $y=x-\zeta>0$.
Differentiating the density (\ref{eq:density-inversion-formula-half}) by
 $\beta$ and doing the same transform as in Lemma \ref{lem:f_alpha(x;alpha,beta)}, we obtain 
\begin{equation*}
f_\beta(x;\alpha,\beta) =
\frac{\tan(\frac{\pi\alpha}{2})}{\pi y}
\left(
\frac{1}{y^\alpha} 
\mbox{Re} \int^\infty_0 e^{i\varsigma}e^{-ite^{i\varsigma}}t^\alpha dt  
-\frac{1}{y} \mbox{Re} \int_0^{\infty}ie^{2i\varsigma}e^{-ite^{i\varsigma}}t dt
+\frac{c}{y^{2\alpha}}+\frac{c}{y^{1+\alpha}} 
\right).
\end{equation*}
The calculation of $f_\beta(x;\alpha,\beta)$ is easy by the equations (\ref{eq:proof-lem3.1-1}) and
 (\ref{eq:proof-lem3.1-2}) in appendix \ref{ap:proof-of-lemma3.1}. As
 $\Delta \to 0$, $\tan\frac{\pi\alpha}{2}=-\frac{\pi}{2}\Delta+o(\Delta^2)$.
Finally we utilize the relation
 $f_\beta(x;\alpha,\beta)=-f_\beta(x;\alpha,-\beta)$ and the first part of the 
 lemma is proved. 
The second part of the lemma is obvious from 
\begin{equation*}
|f_\beta(x;\alpha,\beta)|\le
 \frac{|\tan\frac{\alpha\pi}{2}|}{\pi}\int_0^{\infty}
e^{-t^\alpha}|t^\alpha-t| dt.
\end{equation*}

 \hfill $\Box$ 
\end{pfl}
\section{Information matrix of general stable distributions}
The Fisher information matrix of general stable distributions are derived in this section.
Table \ref{tbl:infomation matrix-limit} gives limiting values of the 
information matrix at $\alpha=2$. To the author's knowledge other than
the diagonal element and
$I_{\mu\sigma}$, these limiting results have not been given in
literature. In following theorem we obtain asymptotic behavior of
information matrix in more detail for the important diverging cases.
\begin{thm}
As $\Delta=2-\alpha \to 0 $, behavior of the Fisher information matrix
 of general stable distributions at $\mu=0$, $\sigma=1$, $\beta \neq
 \pm1$, is given as follows.
\begin{eqnarray*}
\left[
\begin{array}{cccc}
I_{\mu\mu}  & I_{\mu\sigma}    & I_{\mu\alpha}    & I_{\mu\beta} \\
\ast        & I_{\sigma\sigma} & I_{\sigma\alpha} & I_{\sigma\beta} \\
\ast        & \ast             & I_{\alpha\alpha}  & I_{\alpha\beta}  \\
\ast        & \ast             & \ast             & I_{\beta\beta}
\end{array}
\right] =
\left[
\begin{array}{cccc}
0.5+o(1)  & o(1)   & o(1)   & O(\Delta) \\
&&&\\
\ast   & 2.0+o(1)     & -\frac{1}{2}\log\log1/\Delta  & o(\Delta \log\log 1/\Delta) \\
          &              &                {\small \times  (1+o(1))}
	   &  \\
&&& \\
\ast      & \ast             & \frac{1}{4\Delta \log1/\Delta}   & o\left(\frac{1}{\log1/\Delta}\right)  \\
                         &             &  {\small \times(1+o(1))}                     &\\
&&& \\
\ast      & \ast             & \ast             &
 \frac{\Delta}{4(1-\beta^2)\log1/\Delta} \\
&&           &   {\small \times (1+o(1))}   \\
\end{array}
\right].
\end{eqnarray*}
\end{thm}

Note that for $I_{\mu\sigma}$ and $I_{\mu\alpha}$ the rates of convergences
to 0 are not obtained and only the limiting values are
obtained. For $I_{\sigma\alpha}$, $I_{\alpha\alpha}$ and
$I_{\alpha\beta}$ the exact rates of convergences are obtained.
However for symmetric stable distributions since
$I_{\mu\sigma}=I_{\mu\alpha}=0$, the above $3\times3$ matrix presents exact
limiting behavior of the Fisher information matrix as $\Delta \to 0$.
The idea of the proof is based on the proof of Theorem 2 of Nagaev  and  Shkol'nik (1988). 
 
\begin{table}
\caption[percentage]{Limit of information matrix at $\alpha=2$} 
\vspace{-3mm}
\label{tbl:infomation matrix-limit}
\begin{center}
\begin{tabular}{|c|cccc|} \hline
 $I_{\theta\theta}$      & $\mu$ & $\sigma$ & $\alpha$ & $\beta$  \\ \hline
$\mu$    &  0.5  &   0      &    0     &   0      \\
$\sigma$ &   0   &   2.0    & $-\infty$ &   0       \\
$\alpha$ &   0   & $-\infty$ & $\infty$ &   0       \\
$\beta$  &   0   &    0     &   0      &   0   \\  \hline  
\end{tabular}
\end{center}
\end{table}

\begin{pft}
In our proof we use the following notations.
$T$ is a sufficiently large constant such that Theorem
 \ref{thm:density-general-stable-distributions} and Corollaries 
\ref{cor:score-location} and \ref{cor:score-scale} are applicable.
We also denote
\begin{eqnarray*}
&&
x_1(\Delta)=(2-\delta)(\log 1/\Delta)^{1/2}, \quad 
x_2(\Delta)=(2+\delta)(\log 1/\Delta)^{1/2}, \quad
x_3(\Delta)=\exp(\Delta^{-1/2}), 
\end{eqnarray*}
 where $\delta>0$ is an arbitrarily small constant. Further 
we use the notation ``$\mbox{const}$'' for some proper constants. 
We prove $I_{\alpha\alpha}$ in some detail and the proof of the other $I_{\theta\theta}$
 are given in Appendix \ref{ap:pf-information} . \vspace{2mm}\\
{\bf Proof of} $I_{\alpha\alpha}$ :\\
In our proof we utilize Theorem \ref{thm:density-general-stable-distributions} and Lemma \ref{lem:f_alpha(x;alpha,beta)}. 
We divide integral of Fisher information matrix into two parts,
\begin{eqnarray*}
I_{\alpha\alpha} &=&
\int_0^\infty\frac{\{f_\alpha(x+\zeta;\alpha,\beta)\}^2}{f(x+\zeta;\alpha,\beta)} dx
+
\int_0^\infty\frac{\{f_\alpha(x-\zeta;\alpha,-\beta)\}^2}{f(x-\zeta;\alpha,-\beta)}dx
\\
&=& I_{\alpha\alpha}^1+I_{\alpha\alpha}^2. \nonumber
\end{eqnarray*}
This is obtained by the relation $f(x;\alpha,\beta)=f(-x;\alpha,-\beta)$
 for $x-\zeta<0$. 
$I_{\alpha\alpha}^1$ is calculated first. Further we divide integration
$I_{\alpha\alpha}^1$ into five subintegrals,
\begin{equation*}
I_{\alpha\alpha}^1=\sum_{k=1}^5 I_{\alpha\alpha}(k), 
\end{equation*}
where each $I_{\alpha\alpha}(k)$ corresponds to the $k$-th interval of $[0, T)$,\
 $[T,x_1(\Delta))$,\ $[x_1(\Delta),x_2(\Delta))$,\
 $[x_2(\Delta),x_3(\Delta))$ and $[x_3(\Delta),\infty)$. 
 
Clearly 
\begin{equation}
\label{eq:I_{alphaalpha}(1)}
I_{\alpha\alpha}(1) < \infty. 
\end{equation}
For $I_{\alpha\alpha}(2)$ 
\begin{equation*}
f(x+\zeta;\alpha,\beta)=f(x;2)(1+o(1)),\quad
 f_\alpha(x+\zeta;\alpha,\beta)=\mbox{const}\times x^{\Delta-3}.
\end{equation*}
 Then
\begin{equation}
I_{\alpha\alpha}(2)
= \mbox{const}\times
 \int_T^{x_1(\Delta)}x^{2\Delta-6}\exp\left(\frac{x^2}{4}\right) dx
\le \mbox{const}\times
 \frac{1}{\Delta^{1-\delta}(\log 1/\Delta)^{5/2-\Delta}}.
\end{equation}
For $I_{\alpha\alpha}(3)$ 
\begin{equation*}
f(x+\zeta;\beta,\alpha)=
 \{(1+\beta)\Delta x^{\Delta-3}+f(x;2)\}(1+o(1)),\quad
 f_\alpha(x+\zeta;\alpha,\beta)=\mbox{const}\times x^{\Delta-3}.
\end{equation*}
 Then
\begin{equation}
I_{\alpha\alpha}(3)\le\mbox{const}\times \frac{1}{\Delta}
\int_{x_1(\Delta)}^{x_2(\Delta)}x^{\Delta-3}dx 
=
\mbox{const}\times\frac{\delta}{\Delta\log1/\Delta}. 
\end{equation}
For $I_{\alpha\alpha}(4)$ 
\begin{equation*}
f(x+\zeta;\alpha,\beta)=(1+\beta)
 \Delta x^{\Delta-3}(1+o(1)), \quad
 f_\alpha(x+\zeta;\alpha,\beta)=-(1+\beta)x^{\Delta-3}(1+o(1)).
\end{equation*}
 Then
\begin{equation}
I_{\alpha\alpha}(4)=
\frac{1+\beta}{\Delta}\int_{x_2(\Delta)}^{x_3(\Delta)}x^{\Delta-3}
dx (1+o(1)) 
=
\frac{1+\beta}{8\Delta\log1/\Delta}(1+o(1)).
\end{equation}
For $I_{\alpha\alpha}(5)$ 
\begin{equation*}
f(x+\zeta;\alpha,\beta)=(1+\beta)\Delta
 x^{\Delta-3}(1+o(1)),\quad
f_\alpha(x+\zeta;\alpha,\beta)=-(1+\beta+\Delta\log
 x)x^{\Delta-3}(1+o(1)).
\end{equation*}
Then
\begin{eqnarray}
\label{eq:I_{alphaalpha}(5)}
I_{\alpha\alpha}(5)
&=&
\mbox{const}\times \frac{1}{\Delta}\int_{x_3(\Delta)}^{\infty}
x^{\Delta-3}
\{\max(1+\beta,\Delta \log x)\}^2 dx \\
&\le&
\mbox{const}\times \Delta \int_{x_3(\Delta)}^{e^{(1+\beta)/\Delta}}
x^{\Delta-3}(\log x)^2 dx 
+
\mbox{const}\times \frac{1}{\Delta}\int_{e^{(1+\beta)/\Delta}}^\infty
x^{\Delta-3} dx \nonumber \\
&\le&
\mbox{const}\times \frac{1}{\Delta}\int_{x_3(\Delta)}^\infty
x^{\Delta-3} dx \nonumber \\
&=&
O(e^{-2/\Delta^{1/2}}/\Delta)\ \to 0,\ \mbox{as}\ \Delta\ \to 0. \nonumber
\end{eqnarray}
From the formulas (\ref{eq:I_{alphaalpha}(1)})-(\ref{eq:I_{alphaalpha}(5)})
 we obtain
\begin{equation*}
I_{\alpha\alpha}^1=\frac{1+\beta}{8\Delta \log 1/\Delta}(1+o(1)).
\end{equation*}
Setting $\beta \to -\beta$ in $I_{\alpha\alpha}^1$, we obtain
 $I_{\alpha\alpha}^2$. Adding $I_{\alpha\alpha}^1$ and
 $I_{\alpha\alpha}^2$, we prove the assertion.
\\ \hfill $\Box$
\end{pft}

\appendix
\section{Proof of Lemma \ref{lem:1-A}}
As $\Delta \to 0$, $\alpha\varrho= \frac{2}{\pi} \arctan \left(\beta \tan
(\frac{\pi\alpha}{2})\right)$ in $A(\varphi,\alpha,\beta)$ 
 is expanded as
\begin{equation*}
\alpha\varrho=\frac{2}{\pi} \arctan \left(\beta\tan\left(\frac{\pi\alpha}{2}\right)\right) =
-\beta\Delta-\frac{\beta(1-\beta^2)}{3}\left(\frac{\pi}{2}\right)^2\Delta^3+o(\Delta^4).
\end{equation*}
Since the smaller remainder terms other
than $-\beta\Delta$ are not needed in the following argument
we substitute $-\beta\Delta$ into $A(\varphi,\alpha,\beta)$ for
$\alpha \varrho$ in advance. 
Then for convenience we write $A(\varphi,\alpha,\beta)$ of (\ref{integrand-density-A1}) as 
\begin{eqnarray}
\label{integrand-density-A2}
A(\varphi,\alpha,\beta)&=&
\frac{\cos(\frac{\pi}{2}\{(1-\Delta)\varphi-\beta\Delta\})}
{\sin(\frac{\pi}{2}\{(2-\Delta)\varphi-\beta\Delta\})}
\left(
\frac{\cos(\frac{\pi}{2}\varphi) \cos(\frac{\pi}{2}\beta\Delta)}
{\sin(\frac{\pi}{2}\{(2-\Delta)\varphi-\beta\Delta\})}
\right)^{1/(1-\Delta)} \\
&=&
\frac{C(\varphi,\Delta,\beta)}{B(\varphi,\Delta,\beta)}
\left(\frac{D(\varphi)E(\Delta,\beta)}{B(\varphi,\Delta,\beta)}\right)^{1/(1-\Delta)}, \nonumber
\end{eqnarray} 
where 
\begin{eqnarray*}
B(\varphi,\Delta,\beta)&=&
\sin\left(\frac{\pi}{2}\{(2-\Delta)\varphi-\beta\Delta\}\right) , \\
C(\varphi,\Delta,\beta)&=&
\cos\left(\frac{\pi}{2}\{(1-\Delta)\varphi-\beta\Delta\}\right), \\
D(\varphi) &=&
\cos\left(\frac{\pi}{2}\varphi\right), \\
E(\Delta,\beta) &=&
\cos\left(\frac{\pi}{2}\beta\Delta\right).
\end{eqnarray*}
Then substituting $1-\lambda$ for $\varphi$ we obtain
\begin{eqnarray}
B(1-\lambda,\Delta,\beta)&=&\label{eq:expansion-B-pf-lemma1}
 \sin\left(\frac{\pi}{2}(\Delta+2\lambda-\Delta\lambda+\beta\Delta)\right), \\
C(1-\lambda,\Delta,\beta)&=&
\sin\left(\frac{\pi}{2}(\Delta+\lambda-\Delta\lambda+\beta\Delta)\right), \\
D(1-\lambda) &=&
\sin\left(\frac{\pi}{2}\lambda\right),\label{eq:expansion-D-pf-lemma1} \\
E(\Delta,\beta) &=&
\cos\left(\frac{\pi}{2}\beta\Delta\right).
\end{eqnarray}
Note that if $\Delta \to 0$, then $\lambda \to 0$ since $0 \le \lambda \le
\Delta/\epsilon'$ and $\epsilon'>0$ is an arbitrarily
small fixed number (or even when $\epsilon'$ converges to 0 slower than
$\Delta$). Expanding in terms of $\Delta \to 0$, we have
\begin{eqnarray*}
B(1-\lambda,\Delta,\beta)&=&
 \frac{\pi}{2}\Delta\left(1+2\frac{\lambda}{\Delta}-\lambda+\beta\right)+O(\Delta^3) , \\
C(1-\lambda,\Delta,\beta)&=&
 \frac{\pi}{2}\Delta\left(1+\frac{\lambda}{\Delta}-\lambda+\beta\right)+O(\Delta^3), \\
D(1-\lambda) &=&
 \frac{\pi}{2}\Delta \times \frac{\lambda}{\Delta}+O(\Delta^3), \\
E(\Delta,\beta) &=&
1-\frac{1}{2}\left(\frac{\pi}{2}\beta\Delta\right)^2+o(\Delta^3).
\end{eqnarray*} 
Substituting these expansions to $A(\varphi,\Delta,\beta)$ of (\ref{integrand-density-A2}) and 
utilizing $x^{1/(\Delta-1)}=x^{-1}-x^{-1}\log x \Delta+o(\Delta)$,
we prove the desired result.\\
 \hfill $\Box$ 

\section{Proof of Lemma \ref{lem:2-A}}
We prove $A(\varphi_\Delta,\Delta,\beta)$ first.
When $\varphi=\varphi_\Delta=1-\Delta^{1/2-\epsilon}$, we
replace $\lambda$ by $\Delta^{1/2-\epsilon}$ in the formulas 
(\ref{eq:expansion-B-pf-lemma1})-(\ref{eq:expansion-D-pf-lemma1}).
Then
\begin{eqnarray*}
B(\varphi_\Delta,\Delta,\beta)&=&
 \sin\left(\frac{\pi}{2}\{2\Delta^{1/2-\epsilon}+(1+\beta)\Delta-\Delta^{3/2-\epsilon}\}\right), \\
C(\varphi_\Delta,\Delta,\beta)&=&
\sin\left(
     \frac{\pi}{2}\{\Delta^{1/2-\epsilon}+(1+\beta)\Delta-\Delta^{3/2-\epsilon}\}\right), \\
D(\varphi_\Delta) &=&
\sin\left(\frac{\pi}{2}\Delta^{1/2-\epsilon}\right).
\end{eqnarray*}
Note that $E(\Delta,\beta)$ is the same as in Proof of Lemma \ref{lem:1-A} and we omit it. 
As $\Delta \to 0$ we have
\begin{eqnarray}
&&B(\varphi_\Delta,\Delta,\beta) \label{eq:expansion-B-pf-lemma2}\\
&&\quad=
 \pi\Delta^{1/2-\epsilon}\left\{1+\frac{1}{2}(1+\beta)\Delta^{\frac{1}{2}+\epsilon}-\frac{\pi^2}{6}\Delta^{1-2\epsilon}-\frac{1}{2}\Delta 
 -\frac{\pi^2}{4}(1+\beta)\Delta^{3/2-\epsilon}  +o(\Delta^{3/2})\right\}, \nonumber \\
&&C(\varphi_\Delta,\Delta,\beta) \\
&&\quad=
 \frac{\pi}{2}\Delta^{1/2-\epsilon}
\left\{1+(1+\beta)\Delta^{1/2+\epsilon}-
 \frac{1}{6}\frac{\pi^2}{4}\Delta^{1-2\epsilon}-\Delta
-\frac{\pi^2}{8}(1+\beta)\Delta^{3/2-\epsilon}
 +o(\Delta^{3/2}) \right \},\nonumber \\
&&D(\varphi_\Delta)=
\frac{\pi}{2}\Delta^{1/2-\epsilon}\left(1-\frac{1}{6}\frac{\pi^2}{4}\Delta^{1-2\epsilon}
+o(\Delta^{3/2})\right).\label{eq:expansion-D-pf-lemma2}
\end{eqnarray}
Substituting these expansions into $A(\varphi,\Delta,\beta)$ of (\ref{integrand-density-A2}),
we obtain the desired result.

Next we prove assertion concerning $A'(\varphi_\Delta,\Delta,\beta)$.
From the equation (\ref{integrand-density-A1}), we can write
\begin{eqnarray}
&& A'(\varphi,\Delta,\beta)=
A(\varphi,\Delta,\beta) \label{first-derivative-integrand-density}  \\
&& \times \frac{\pi}{2}
\left\{
\frac{(2-\Delta)^2}{\Delta-1}\frac{\cos(\frac{\pi}{2}\{(2-\Delta)\varphi-\beta\Delta\})}{\sin(\frac{\pi}{2}\{(2-\Delta)\varphi-\beta\Delta\})}
-(1-\Delta)\frac{\sin(\frac{\pi}{2}\{(1-\Delta)\varphi-\beta\Delta\})}{\cos(\frac{\pi}{2}\{(1-\Delta)\varphi-\beta\Delta\})}
-\frac{1}{1-\Delta}\frac{\sin(\frac{\pi}{2}\varphi)}{\cos(\frac{\pi}{2}\varphi)}
\right\} \nonumber \\
&&=A(\varphi,\Delta,\beta) \times
\frac{\pi}{2}\left\{
\frac{(2-\Delta)^2}{\Delta-1}\frac{F(\varphi,\Delta,\beta)}{B(\varphi,\Delta,\beta)}-(1-\Delta)\frac{G(\varphi,\Delta,\beta)}{C(\varphi,\Delta,\beta)}
-\frac{1}{1-\Delta}\frac{H(\varphi)}{D(\varphi)}
\right\}, \nonumber
\end{eqnarray}
where
\begin{eqnarray*}
F(\varphi,\Delta,\beta)&=&
\cos\left(
     \frac{\pi}{2}\{(2-\Delta)\varphi-\beta\Delta\}\right), \\
G(\varphi,\Delta,\beta)&=&
\sin\left(\frac{\pi}{2}\{(1-\Delta)\varphi-\beta\Delta\}\right), \\
H(\varphi) &=&
\sin\left(\frac{\pi}{2}\varphi\right).
\end{eqnarray*}
Then substituting $\varphi_\Delta$ for $\varphi$, we obtain
\begin{eqnarray*}
F(\varphi_\Delta,\Delta,\beta)&=&
-\cos\left(
     \frac{\pi}{2}\{2\Delta^{1/2-\epsilon}+(1+\beta)\Delta-\Delta^{3/2-\epsilon}\}\right), \\
G(\varphi_\Delta,\Delta,\beta)&=&
\cos\left(
     \frac{\pi}{2}\{\Delta^{1/2-\epsilon}+(1+\beta)\Delta-\Delta^{3/2-\epsilon}\}\right), \\
H(\varphi_\Delta) &=&
\cos\left(\frac{\pi}{2}\Delta^{1/2-\epsilon}\right).
\end{eqnarray*}
As $\Delta \to 0$, 
\begin{eqnarray}
F(\varphi_\Delta,\Delta,\beta)&=&\label{eq:expansion-F-pf-lemma2}
 -1+\frac{1}{2}\pi^2\Delta^{1-2\epsilon}+\frac{1}{2}\pi^2(1+\beta)\Delta^{3/2-\epsilon}+o(\Delta^{3/2}) , \\
G(\varphi_\Delta,\Delta,\beta)&=&
 1-\frac{\pi^2}{8}\Delta^{1-2\epsilon}-\frac{\pi^2}{4}(1+\beta)\Delta^{3/2-\epsilon}+o(\Delta^{3/2}) , \\
H(\varphi_\Delta) &=&
1-\frac{\pi^2}{8}\Delta^{1-2\epsilon}+o(\Delta^{3/2}).\label{eq:expansion-H-pf-lemma2}
\end{eqnarray}
Combining
 (\ref{eq:expansion-B-pf-lemma2})-(\ref{eq:expansion-D-pf-lemma2}) ,
 (\ref{eq:expansion-F-pf-lemma2})-(\ref{eq:expansion-H-pf-lemma2}),
 $A(\varphi,\Delta,\beta)$ and $A'(\varphi,\alpha,\beta)$ of
(\ref{first-derivative-integrand-density}), we prove the lemma.\\
\hfill $\Box$ 

\section{Proof of Lemma \ref{lem:3-A}}
We further expand $B(\varphi,\Delta,\beta)$, $C(\varphi,\Delta,\beta)$ and $D(\varphi,\Delta)$ in
(\ref{eq:expansion-B-pf-lemma1})-(\ref{eq:expansion-D-pf-lemma1}) with
respect to $\lambda$. As
$\Delta \to 0$ and $\lambda \to 0$ such that $\Delta^{1/2-\epsilon}
\le \lambda$,
\begin{eqnarray}
 &&\hspace{-3mm}B(1-\lambda,\Delta,\beta)\label{eq:expansion-B-pf-lemma3}  \\
&&=\pi\lambda\left\{
1-\frac{1}{6}\left(\pi\lambda\right)^2
+\frac{1}{120}(\pi\lambda)^4
+\frac{\Delta(1-\lambda+\beta)}{2\lambda}
-\frac{\pi^2\lambda\Delta}{4}(1-\lambda+\beta)
+o(\lambda^4)
\right\}, \nonumber \\
&&\hspace{-3mm}C(1-\lambda,\Delta,\beta) \\
&&= \frac{\pi}{2}\lambda\left\{
1-\frac{1}{6}\left(\frac{\pi}{2}\lambda\right)^2
+\frac{1}{120}\left(\frac{\pi}{2}\lambda\right)^4
+\frac{\Delta(1-\lambda+\beta)}{\lambda}
-\frac{\pi^2\lambda\Delta}{8}(1-\lambda+\beta)
+o(\lambda^4)
\right\}, \nonumber \\
&&\hspace{-3mm}D(1-\lambda,\Delta) \label{eq:expansion-D-pf-lemma3}\\
&&=
\frac{\pi}{2}\lambda\left\{1-\frac{1}{6}\left(\frac{\pi}{2}\lambda\right)^2
+\frac{1}{120}\left(\frac{\pi}{2}\lambda\right)^4
+o(\lambda^5)\right\}. \nonumber
\end{eqnarray}
Substituting these equations to $A(\varphi,\Delta,\beta)$, we obtain the
result.

Concerning $A'(\varphi,\alpha,\beta)$ we use formulas
\begin{eqnarray*}
F(1-\lambda,\Delta,\beta)&=& 
-\cos\left(
     \frac{\pi}{2}\{2\lambda+(1-\lambda+\beta)\Delta \}\right), \\
G(1-\lambda,\Delta,\beta)&=&
\cos\left(
     \frac{\pi}{2}\{\lambda+(1-\lambda+\beta)\Delta\}\right), \\
H(1-\lambda) &=&
\cos\left(\frac{\pi}{2}\lambda\right). 
\end{eqnarray*}
As
$\Delta \to 0$ and $\lambda \to 0$ under $\Delta^{1/2-\epsilon}
\le \lambda$,
\begin{eqnarray}
F(1-\lambda,\Delta,\beta)&=& \label{eq:expansion-F-pf-lemma3}
-1+\frac{1}{2}(\pi\lambda)^2+o(\lambda^3), \\
G(1-\lambda,\Delta,\beta)&=&
1-\frac{1}{2}\left(\frac{\pi}{2}\lambda\right)^2+o(\lambda^3), \\
H(1-\lambda) &=&
1-\frac{1}{2}\left(\frac{\pi}{2}\lambda\right)^2+o(\lambda^3).\label{eq:expansion-H-pf-lemma3}
\end{eqnarray}
Substituting (\ref{eq:expansion-B-pf-lemma3})-(\ref{eq:expansion-D-pf-lemma3}) and 
(\ref{eq:expansion-F-pf-lemma3})-(\ref{eq:expansion-H-pf-lemma3}) to
$A'(\varphi,\alpha,\beta)$ of
(\ref{first-derivative-integrand-density}), we obtain the desired
assertion.\\
From the equation (\ref{integrand-density-A1}) or $A'(\varphi,\alpha,\beta)$ of
(\ref{first-derivative-integrand-density}), we can write
\begin{eqnarray}
&& A''(\varphi,\Delta,\beta)  \label{second-derivative-integrand-density} \\
&&=A'(\varphi,\Delta,\beta) \times
\frac{\pi}{2}\left\{
\frac{(2-\Delta)^2}{\Delta-1}\frac{F(\varphi,\Delta,\beta)}{B(\varphi,\Delta,\beta)}-(1-\Delta)\frac{G(\varphi,\Delta,\beta)}{C(\varphi,\Delta,\beta)}
-\frac{1}{1-\Delta}\frac{H(\varphi)}{D(\varphi)}
\right\} \nonumber\\
&&+A(\varphi)\frac{\pi^2}{4(1-\Delta)}
\left\{
\frac{(2-\Delta)^3}{\left(\sin(\frac{\pi}{2}\{(2-\Delta)\varphi-\beta\Delta\})\right)^2}
-\frac{(1-\Delta)^3}{\left(\cos(\frac{\pi}{2}\{(1-\Delta)\varphi-\beta\Delta\})\right)^2}
-\frac{1}{\left(\cos( \frac{\pi}{2}\varphi)\right)^2 \}}
\right\} \nonumber \\
&&=A'(\varphi,\Delta,\beta) \times
\frac{\pi}{2}\left\{
\frac{(2-\Delta)^2}{\Delta-1}\frac{F(\varphi,\Delta,\beta)}{B(\varphi,\Delta,\beta)}-(1-\Delta)\frac{G(\varphi,\Delta,\beta)}{C(\varphi,\Delta,\beta)}
-\frac{1}{1-\Delta}\frac{H(\varphi)}{D(\varphi)}
\right\} \nonumber\\
&&+A(\varphi)\frac{\pi^2}{4(1-\Delta)}
\left\{
\frac{(2-\Delta)^3}{B^2(\varphi,\Delta,\beta)}
-\frac{(1-\Delta)^3}{C^2(\varphi,\Delta,\beta)}
-\frac{1}{D^2(\varphi)}
\right\}. \nonumber
\end{eqnarray}
Combining (\ref{eq:expansion-B-pf-lemma3})-(\ref{eq:expansion-D-pf-lemma3}), 
(\ref{eq:expansion-F-pf-lemma3})-(\ref{eq:expansion-H-pf-lemma3}), 
$A'(\varphi,\alpha,\beta)$ of (\ref{first-derivative-integrand-density}) and $A''(\varphi,\alpha,\beta)$ of
(\ref{second-derivative-integrand-density}), we obtain the desired
assertion. \\
\hfill $\Box$

\section{Proof of Lemma \ref{lem:4-A}}
We only need to prove for $x-\zeta>0$, because for $x-\zeta<0$
$f(x;\alpha,\beta)=f(-x;\alpha,-\beta)$ and $\varrho^{\ast}$ in 
(\ref{integrand-density-A1}) does not change, namely
\begin{equation*}
\varrho^{\ast}=\sgn(x-\zeta)\frac{\pi\alpha}{2}\arctan\left(\beta\tan\left(\frac{\pi\alpha}{2}\right)\right)
=\sgn(-x+\zeta)\frac{\pi\alpha}{2}\arctan\left(-\beta\tan\left(\frac{\pi\alpha}{2}\right)\right).
\end{equation*}
 We use notation $\beta_B \in (-1,1)$ which corresponds to the parameter
 $\beta \in (-1,1)$ for $(B)$ representation of characteristic
 functions. The reason is that the proof in terms of (B) representation is simpler since $\alpha\varrho=-\beta_B\Delta$
(see p.74 of Zolotarev 1986 ; the notation $\varrho$ in our paper
 corresponds to  $\theta$ in Zolotarev) and we do not need to expand $\alpha\varrho$ with respect
 to $\beta$ as in the proof of Lemma \ref{lem:3-A}. Moreover our purpose is only to determine the sign of $A'(\varphi)$. 
We divide the integral interval as $(-\varrho,1-\Delta/\epsilon']$,
 $[1-\Delta/\epsilon',1)$ and prove the lemma for each interval.
Except for the interval $(-\varrho,0)$ we follow the line of the proof of Lemma 4 of Nagaev and Shkol'nik (1988).
For $(1-\Delta/\epsilon',1)$ the proof is obvious from Lemma \ref{lem:1-A}.
For $(-\varrho,1-\Delta/\epsilon']$ we have to consider the sign of $\varrho$.
Note that
\begin{equation*}
\sgn(\varrho)=-\sgn(\beta_B).
\end{equation*}
First for $0 \le -\varrho \le \varphi \le 1-\Delta/\epsilon'\ \Leftrightarrow\ \Delta/\epsilon' \le \lambda
\le 1+\varrho$, we get 
\begin{equation*} 
A'(\varphi) =\frac{-1}{1-\Delta} \frac{\pi}{2} A(\varphi)A_1(\varphi),
\end{equation*}
where
\begin{eqnarray*}
A_1(1-\lambda)&=&
-(2-\Delta)^2
\cot\left(\frac{\pi}{2}(\Delta+2\lambda-\Delta\lambda+\beta_B\Delta)\right)\\
&&+(1-\Delta)^2\cot\left(\frac{\pi}{2}(\Delta+\lambda-\Delta\lambda+\beta_B\Delta)\right)+\cot\left(\frac{\pi}{2}
\lambda\right). 
\end{eqnarray*}
From the formula (4.3.70) in [4] we have 
\begin{equation*}
A_1(1-\lambda)=
\frac{2(1+\beta_B)^2\Delta^2}{\pi\lambda(\Delta+2\lambda-\Delta\lambda+\beta_B\Delta)(\Delta+\lambda-\Delta\lambda+\beta_B\Delta)}
+\sum_{k=0}^{\infty}d_k  \gamma_k(\lambda),
\end{equation*}
where $d_k$'s are positive numbers and
\begin{equation*}
\gamma_k(\lambda)
=\left(\frac{\pi}{2} \lambda \right)^{2k+1}
\left[
(2-\Delta)^2\{2-\Delta+(1+\beta_B)\Delta/\lambda \}
-(1-\Delta)^2\{1-\Delta+(1+\beta_B)\Delta/\lambda \}
\right].
\end{equation*}
Since $\Delta/\lambda \le \epsilon'$ for arbitrarily small $\epsilon'>0$,
$A_1(\varphi)>0$ is easily confirmed. 

Secondly we investigate the case $-\varrho < \varphi
\le 0 \ \Leftrightarrow \alpha \varphi=\gamma\Delta$ where $\gamma \in
(\beta_B,0)$ for $\beta_B<0$.  Then
\begin{eqnarray*}
A_1(\varphi) 
&=&
(2-\Delta)^2\cot\left( \frac{\pi}{2}\alpha(\varphi+\varrho)\right)
+(1-\Delta)^2\tan\left(\frac{\pi}{2}\{\alpha\varrho+(\alpha-1)\varphi \}\right)
+\tan\left(\frac{\pi}{2}\varrho\right)  \\
&=&
(2-\Delta)^2\cot\left(\frac{\pi}{2}(\gamma-\beta_B)\Delta\right)
+(1-\Delta)^2\tan\left(\frac{\pi}{2}\left(\gamma-\beta_B-\frac{\gamma}{2-\Delta}\right)\Delta\right) 
+\tan\left(\frac{\pi}{2}\frac{\gamma\Delta}{2-\Delta}\right).
\end{eqnarray*}
As $\Delta \to 0$,
\begin{equation*}
A_1(\varphi)=
\frac{(2-\Delta)^2}{\frac{\pi}{2}(\gamma-\beta_B)\Delta}
+O (\Delta).
\end{equation*} 
We easily find that $A_1(\varphi)>0$ for sufficiently small $\Delta$.
\\
 \hfill $\Box$ 
\section{Calculations of $f_\alpha$}
\label{ap:proof-of-lemma3.1}
For all $y \ge y_0=x_0-\zeta$, we have
\begin{equation*}
\exp\left(it^{\alpha}y^{-\alpha}\left(1-i\beta\tan\left(\frac{\pi\alpha}{2}\right)\right)\right)
=1+c_0\ t^{\alpha}y^{-\alpha}\left(1-i\beta\tan\left(\frac{\alpha\pi}{2}\right)\right),
\end{equation*}
where $c_0 \in \mbox{\bf C}$ is some complex constant.
Utilizing this expansion in $f_\alpha(x;\alpha,\beta)$ we obtain
\begin{eqnarray*}
K_1 &=&
\frac{1}{\pi y}\mbox{Re}\ e^{i\varsigma}
\int^{\infty}_{0}\exp\left(-ite^{i\varsigma}+it^{\alpha}y^{-\alpha}\left(1-i\beta\tan\left(\frac{\pi\alpha}{2}\right)\right)\right)
\nonumber \\
&&\times  it^{\alpha}y^{-\alpha}(\log t+i\varsigma+\log
 y)\left(1-i\beta\tan\left(\frac{\pi\alpha}{2}\right)\right) dt  \\
&=&
\frac{1}{\pi y^{1+\alpha}}\mbox{Re}\ e^{i\varsigma}\int_0^\infty 
e^{-ite^{i\varsigma}}\left\{1+c_0\ t^{\alpha}y^{-\alpha}\left(1-i\beta\tan\left(\frac{\alpha\pi}{2}\right)\right)\right\}\\
&& \times it^\alpha(\log t+i\varsigma+\log y)\left(1-i\beta\tan\left(\frac{\pi\alpha}{2}\right)\right)dt \\
&=&
\frac{1}{\pi y^{1+\alpha}} \mbox{Re}\ e^{i\varsigma}\int_0^\infty 
e^{-ite^{i\varsigma}}it^{\alpha}(\log t+i\varsigma+\log
y)\left(1-i\beta\tan\left(\frac{\alpha\pi}{2}\right)\right) dt \\
&&+\frac{1}{\pi y^{1+2\alpha}} \mbox{Re}\ e^
{i\varsigma}\int_0^\infty
 e^{-ite^{i\varsigma}}
ic_0\,t^{2\alpha}\left(1-\beta\tan\left(\frac{\alpha\pi}{2}\right)\right)^2(\log t+
i\varsigma+\log y)dt \\
&=&
\frac{1}{\pi y^{1+\alpha}}\mbox{Re}\ e^{i\varsigma}\int_0^\infty 
e^{-ite^{i\varsigma}}it^\alpha (\log t+i\varsigma+\log z)dt +
\frac{\Delta(c+c\log y)}{y^{1+\alpha}}+\frac{c+c\log y}{y^{1+2\alpha}}.
\end{eqnarray*}
Note that As
 $\Delta \to 0$,
 $\tan\frac{\pi\alpha}{2}=-\frac{\pi}{2}\Delta+o(\Delta^2)$. Then we
 calculate the integrations in $K_1$ above.
Write $\lambda_1=\cos\varsigma$ and $\lambda_2=\sin
\varsigma$. Here at $\alpha=2$, $\lambda_1=1/\sqrt{2}$ and $\lambda_2=-1/\sqrt{2}$. 
From the formulas (3.5), (4.40) in part I and (3.7), (4.17) in part
I\hspace{-.1em}I of Oberhettinger (1990), we have
\begin{eqnarray*}
\mbox{Re}\ e^{i\varsigma}\int_0^\infty 
e^{-ite^{i\varsigma}}it^\alpha \log t dt 
&=&
\int_0^{\infty}e^{t\lambda_2}t^2\log
t\{\lambda_1\sin(t\lambda_1)-\lambda_2\cos(t\lambda_1)\}dt+\Delta c \\
&=&
4\lambda_1\lambda_2\arctan\left(-\frac{\lambda_1}{\lambda_2}\right)+8\lambda_1\lambda_2(\lambda_2^2-\lambda_1^2)\left(\frac{3}{2}-\mbox{EU}\right)
+\Delta c \\
&=&-\frac{\pi}{2}+\Delta c,
\end{eqnarray*}
\begin{eqnarray}
\label{eq:proof-lem3.1-1}
 -\mbox{Re}\ \varsigma e^{i\varsigma}\int_0^\infty 
e^{-ite^{i\varsigma}}t^\alpha dt 
&=&
-\int_0^{\infty}e^{t\lambda_2}t^2
\{\lambda_2\sin(t\lambda_1)+\lambda_1\cos(t\lambda_1)\}dt+\Delta c \\
&=&
-4\varsigma \lambda_1\lambda_2+\Delta c \nonumber \\
&=&
-\frac{\pi}{2}+\Delta c, \nonumber
\end{eqnarray}
\begin{eqnarray*}
\mbox{Re}\ e^{i\varsigma}\int_0^\infty 
e^{-ite^{i\varsigma}}it^\alpha dt 
&=&
\int_0^{\infty}e^{t\lambda_2}t^2
\{\lambda_1\sin(t\lambda_1)-\lambda_2\cos(t\lambda_1)\}dt+\Delta c \\
&=&
-2\lambda_1^4+2\lambda_2^4+\Delta c\\
&=&\Delta c,
\end{eqnarray*}
where $\mbox{EU}$ means Euler's constant.
Substituting these equations into $K_1$, we obtain
\begin{equation}
\label{eq:K1}
K_1=-\frac{1}{y^{1+\alpha}}\left\{1+\Delta (c +c \log y) + \frac{c+c\log y}{y^{\alpha}}\right\}.
\end{equation}
For $K_2$ we do the similar calculations.
\begin{eqnarray*} 
K_2 &=&
-\frac{1}{\pi y}\mbox{Re}\ e^{i\varsigma}
\int^{\infty}_{0}\exp\left(-ite^{i\varsigma}+it^{\alpha}y^{-\alpha}\left(1-i\beta\tan\left(\frac{\pi\alpha}{2}\right)\right)\right)
\nonumber \\
&&\times 
(it^{\alpha}y^{-\alpha}+ty^{-1}e^{i\varsigma})\frac{i\beta\pi}{2}/\left(\cos\frac{\pi\alpha}{2}\right)^2
dt  \\
&=&
\frac{\beta}{2y (\cos(\frac{\pi\alpha}{2}))^2}
\mbox{Re}\ e^{i\varsigma}\int_0^\infty 
e^{-ite^{i\varsigma}}\left\{1+ct^{\alpha}y^{-\alpha}\left(1-i\beta\tan\left(\frac{\pi\alpha}{2}\right)\right)\right\}
(t^\alpha y^{-\alpha}-ity^{-1}e^{i\varsigma})dt \\
&=&
\frac{\beta}{2y (\cos(\frac{\pi\alpha}{2}))^2}
\left\{
\frac{1}{y^\alpha}\mbox{Re}\ e^{i\varsigma}\int_0^\infty 
e^{-ite^{i\varsigma}}t^\alpha dt 
-
\frac{1}{y}\mbox{Re}\ e^{2i\varsigma}\int_0^\infty 
e^{-ite^{i\varsigma}}it dt
+
\frac{c}{y^{2\alpha}}
+
\frac{c}{y^{1+\alpha}} 
\right\}.
\end{eqnarray*}
Substituting the following two equations into $K_2$, we obtain $K_2$. The first equation is
obtained from formula (\ref{eq:proof-lem3.1-1}) in the calculation of
$K_1$. The second equation is obtained from the 
formulas (3.5) in part I and (3.7) in part I\hspace{-.1em}I of Oberhettinger (1990).
\begin{eqnarray*}
\mbox{Re}\ e^{i\varsigma}\int_0^\infty 
e^{-ite^{i\varsigma}} t^\alpha dt 
&=&
-2+\Delta c,
\end{eqnarray*}
\begin{eqnarray}
\label{eq:proof-lem3.1-2}
 -\mbox{Re}\ ie^{2i\varsigma}\int_0^\infty 
e^{-ite^{i\varsigma}}t dt 
&=&
\int_0^{\infty}e^{t\lambda_2}t
\{2\lambda_1\lambda_2\cos(t\lambda_1)+(\lambda_2^2-\lambda_1^2)\sin(t\lambda_1)\}
 dt \\
&=&
2\lambda_1\lambda_2(\lambda_2^2-\lambda_1^2)+(\lambda_2^2-\lambda_1^2)(-2\lambda_1\lambda_2)
\nonumber \\
&=&
0. \nonumber
\end{eqnarray}
Hence
\begin{equation}
\label{eq:K2}
K_2=\frac{\beta}{y^{1+\alpha}}
\left\{
-1+\Delta c+\frac{c}{y^\alpha}+\frac{c}{y}
\right\}.
\end{equation}

\section{Proof of $I_{\theta\theta}$ other than $I_{\alpha\alpha}$}
\label{ap:pf-information}
Here we give proofs for
$I_{\beta\beta}$,$I_{\alpha\beta}$,$I_{\sigma\alpha}$,$I_{\sigma\beta}$,
$I_{\mu\beta}$ and $I_{\mu\alpha}$. In the following we assume
$T>x_0$ where $x_0$ is defined in Theorem
\ref{thm:density-general-stable-distributions} and Corollaries
\ref{cor:score-location} and \ref{cor:score-scale}.\vspace{2mm} \\
{\bf Proof of} $I_{\beta\beta}$ :\\
Here we utilize Theorem \ref{thm:density-general-stable-distributions}
 and Lemma \ref{lem:f_beta(x;alpha,beta)}.
However the proof of $I_{\beta\beta}$ is essentially the same as the proof of $I_{\alpha\alpha}$ because 
\begin{equation*}
f_\beta(x;\alpha,\beta)=-\frac{\Delta}{1+\beta}\sgn(x-\zeta)f_\alpha(x;\alpha,\beta)(1+o(1))
\end{equation*}
for sufficiently large $|x-\zeta|$. Therefore we omit the detailed
 derivation of $I_{\beta\beta}$. \vspace{2mm} \\
{\bf Proof of} $I_{\alpha\beta}$ : \\
Here we utilize Theorem \ref{thm:density-general-stable-distributions},
 Lemma \ref{lem:f_alpha(x;alpha,beta)}
 and Lemma \ref{lem:f_beta(x;alpha,beta)}.
From the same reason as $I_{\beta\beta}$ we omit the proof. \vspace{2mm}\\
{\bf Proof of} $I_{\sigma\alpha}$ :\\
We make use of Lemma \ref{lem:f_alpha(x;alpha,beta)}
 and Corollary \ref{cor:score-scale}.
We divide integration $I_{\sigma\alpha}$ into two subintegrals,
\begin{eqnarray*}
I_{\sigma\alpha} &=&
\int_0^\infty\frac{f_\sigma(x+\zeta;\alpha,\beta)f_\alpha(x+\zeta;\alpha,\beta)}{f(x+\zeta;\alpha,\beta)} dx
+
\int_0^\infty\frac{f_\sigma(x-\zeta;\alpha,-\beta)f_\alpha(x-\zeta;\alpha,-\beta)}{f(x-\zeta;\alpha,-\beta)}dx \\
&=& I_{\sigma\alpha}^1+I_{\sigma\alpha}^2. \nonumber
\end{eqnarray*}
$I_{\sigma\alpha}^1$ is calculated first. Further we divide integration
$I_{\sigma\alpha}^1$ into four subintegrals,
\begin{equation*}
I_{\sigma\alpha}^1=\sum_{k=1}^4 I_{\sigma\alpha}(k), 
\end{equation*}
where each $I_{\sigma\alpha}(k)$ corresponds to the integration of
$I_{\sigma\alpha}^1$ for the $k$-th interval of $[0, T)$,\
 $[T,x_1(\Delta))$,\ $[x_1(\Delta),x_2(\Delta))$ and 
 $[x_2(\Delta),\infty)$. 
For $I_{\sigma\alpha}(1)$ clearly
\begin{equation}
\label{eq:I_sigmaalpha(1)}
I_{\sigma\alpha}(1)<\infty.
\end{equation}
For $I_{\sigma\alpha}(2)$
\begin{equation*}
\frac{f_\sigma(x+\zeta,\alpha,\beta)}{f(x+\zeta,\alpha,\beta)}=\frac{x^2}{2}(1+o(1)),
 \quad f(x+\zeta,\alpha,\beta)=-(1+\beta)x^{\Delta-3}(1+o(1)).
\end{equation*}
 Then
\begin{eqnarray}
\label{eq:I_sigmaalpha(2)}
I_{\sigma\alpha}(2)&=&-(1+\beta)\int_T^{x(\Delta)}\frac{x^{\Delta-1}}{2}
 dx (1+o(1))
 \\
&=&
-\frac{1+\beta}{4}\log\log 1/\Delta(1+o(1)). \nonumber
\end{eqnarray}
For $I_{\sigma\alpha}(3)$ 
\begin{equation*}
\left|\frac{f_\sigma(x+\zeta,\alpha,\beta)
 }{f(x+\zeta,\alpha,\beta)}\right|\le 
\mbox{const}\times x^2,\quad
 |f_\alpha(x+\zeta;\alpha,\beta)|=(1+\beta)x^{\Delta-3}(1+o(1)),
\end{equation*}
and
\begin{equation*}
\left|\frac{f_\sigma(x+\zeta;\alpha,\beta)f_\alpha(x+\zeta;\alpha,\beta)}{f(x+\zeta;\alpha,\beta)}\right|
=\mbox{const}\times x^{\Delta-1}(1+o(1)).
\end{equation*}
It is easy to see
\begin{equation}
\label{eq:I_sigmaalpha(3)}
I_{\sigma\alpha}(3)=\mbox{const} \times \int_{x_1(\Delta)}^{x_2(\Delta)} x^{\Delta-1} dx = O(\delta).
\end{equation}
For $I_{\sigma\alpha}(4)$
\begin{equation*}
\left|\frac{f_\sigma(x+\zeta;\alpha,\beta)}{f(x+\zeta;\alpha,\beta)}\right|
\le 2(1+o(1)), \quad |f_\alpha(x+\zeta;\alpha,\beta)|=(1+\beta)x^{\Delta-3}(1+o(1)).
\end{equation*}
Thus we easily show
\begin{equation}
\label{eq:I_sigmaalpha(4)}
I_{\sigma\alpha}(4) \to 0, \quad \mbox{as}\ \Delta \to 0.
\end{equation}
From (\ref{eq:I_sigmaalpha(1)}), (\ref{eq:I_sigmaalpha(2)}), (\ref{eq:I_sigmaalpha(3)}) and (\ref{eq:I_sigmaalpha(4)}), we obtain
\begin{equation*}
I_{\sigma\alpha}^1=-\frac{1}{4}(1+\beta)\log\log 1/\Delta (1+o(1)).
\end{equation*}
Setting $\beta \to -\beta$ in $I_{\sigma\alpha}^1$, we obtain $I_{\sigma\alpha}^2$.
 Adding $I_{\sigma\alpha}^1$ and $I_{\sigma\alpha}^2$, we prove the assertion.
\vspace{2mm} \\
{\bf Proof of} $I_{\sigma\beta}$ : \\
Here we utilize Lemma \ref{lem:f_beta(x;alpha,beta)} and Corollary \ref{cor:score-scale}.
From the same reason as $I_{\beta\beta}$ we
 omit the proof. Note that the result can be expected from behavior of
 $\Delta \times I_{\sigma\alpha}$ and the asymmetry of
 $f_\beta(x;\alpha,\beta)$ around $\zeta$, which we see in Lemma \ref{lem:f_beta(x;alpha,beta)}.  \vspace{2mm} \\
{\bf Proof of} $I_{\mu\beta}$ : \\
Here we utilize Lemma \ref{lem:f_beta(x;alpha,beta)} and Corollary \ref{cor:score-location}.
We divide integration $I_{\mu\beta}$ into two subintegrals.
 \begin{eqnarray*}
I_{\mu\beta} &=&
\int_0^\infty\frac{f_\mu(x+\zeta;\alpha,\beta)f_\beta(x+\zeta;\alpha,\beta)}{f(x+\zeta;\alpha,\beta)} dx
+
\int_0^\infty\frac{f_\mu(x-\zeta;\alpha,-\beta)f_\beta(x-\zeta;\alpha,-\beta)}{f(x-\zeta;\alpha,-\beta)}dx \\
&=& I_{\mu\beta}^1+I_{\mu\beta}^2. \nonumber
\end{eqnarray*}
$I_{\mu\beta}^1$ is calculated first. Further we divide integration
$I_{\mu\beta}^1$ into two subintegrals,
\begin{equation*}
I_{\mu\beta}^1=\sum_{k=1}^2 I_{\mu\beta}(k), 
\end{equation*}
where each $I_{\mu\beta}(k)$ corresponds to the integration of
$I_{\mu\beta}^1$ for the $k$-th interval of $[0, T)$ and
 $[T,\infty)$.
For $I_{\mu\beta}(1)$
\begin{equation*}
f_\beta(x+\zeta;\alpha,\beta)=O(\Delta),\quad
\left|\frac{f_\mu(x+\zeta;\alpha,\beta)}{f(x+\zeta;\alpha,\beta)}\right| \le \mbox{const}.
\end{equation*}
Then $I_{\mu\beta}(1)=O(\Delta)$. 
For $I_{\mu\beta}(2)$ it follows that  
\begin{equation*}
\left|\frac{f_\mu(x+\zeta;\alpha,\beta)}{f(x+\zeta;\alpha,\beta)}\right| \le
 \mbox{const}\times x,\quad
 f_\beta(x+\zeta;\alpha,\beta)=\mbox{const} \times \Delta x^{\Delta-3},
\end{equation*}
uniformly in $x\in[T,\infty)$. 
Thus we find $I_{\mu\beta}(2)=O(\Delta)$ easily. Setting $\beta \to -\beta$ in $I_{\mu\beta}^1$, we obtain $I_{\mu\beta}^2$.
 Adding $I_{\mu\beta}^1$ and $I_{\mu\beta}^2$, we obtain the desired
 result.
\vspace{2mm} \\
{\bf Proof of} $I_{\mu\alpha}$ :\\
Here we utilize Lemma \ref{lem:f_alpha(x;alpha,beta)} and Corollary \ref{cor:score-location}.
We divide integration $I_{\sigma\alpha}$ into two subintegrals,
\begin{eqnarray*}x
I_{\mu\alpha} &=&
\int_0^\infty\frac{f_\mu(x+\zeta;\alpha,\beta)f_\alpha(x+\zeta;\alpha,\beta)}{f(x+\zeta;\alpha,\beta)} dx
-\int_0^\infty\frac{f_\mu(x-\zeta;\alpha,-\beta)f_\alpha(x-\zeta;\alpha,-\beta)}{f(x-\zeta;\alpha,-\beta)}dx \\
&=& I_{\mu\alpha}^1+I_{\mu\alpha}^2. \nonumber
\end{eqnarray*}
$I_{\mu\alpha}^1$ is calculated first.
From Corollary \ref{cor:score-location} and Lemma \ref{lem:f_beta(x;alpha,beta)} we can find
 an integrable function $g_1(x)$ such that for $x \in (0,\infty)$ and
 all sufficiently small $\Delta$ 
\begin{equation*}
\left|\frac{f_\mu(x+\zeta;\alpha,\beta)f_\alpha(x+\zeta;\alpha,\beta)}{f(x;\alpha,\beta)}
\right| \le g_1(x).
\end{equation*}
Hence we can apply the dominated convergence theorem to
 $I_{\mu\alpha}^1$. Then
\begin{equation*}
\lim_{\Delta \to 0}I_{\mu\alpha}^1=
\int_0^\infty\frac{f_\mu(x;2)f_\alpha(x;2)}{f(x;2)} dx.
\end{equation*}
We apply the same arguments as $I_{\mu\alpha}^1$ to $I_{\mu\alpha}^2$
 and obtain
\begin{equation*}
\lim_{\Delta \to 0}I_{\mu\alpha}^2=-
\int_0^\infty\frac{f_\mu(x;2)f_\alpha(x;2)}{f(x;2)} dx.
\end{equation*}
 Adding $I_{\mu\alpha}^1$ and $I_{\mu\alpha}^2$, we prove our assertion.
\\
 \hfill $\Box$

\end{document}